\documentclass[11pt,twoside]{amsart}

\usepackage{amssymb}
\usepackage{amsxtra}
\usepackage{graphicx}
\usepackage{epsfig}
\usepackage{dcpic,pictexwd}

\theoremstyle{plain}

\newtheorem{thm}{Theorem}
\newtheorem{theorem}{Theorem}[section]
\newtheorem{prop}[theorem]{Proposition}

\newtheorem{lemma}{Lemma}
\newtheorem{cor}{Corollary}

\newcommand{\R}{\mathbb{R}}

\newcommand{\Z}{\mathbb{Z}}

\newcommand{\F}{\mathbb{F}}

\newcommand{\Si}{\Sigma}

\newcommand{\ben}{\begin{enumerate}}
\newcommand{\een}{\end{enumerate}}

\newcommand{\SP}[1]{\mathfrak{#1}}

\newcommand{\ra}{\rightarrow}
\newcommand{\lra}{\longrightarrow}
\newcommand{\mir}[1]{\overline{#1}}

\newcommand{\Oz}{P.\ Ozsv{\'a}th\,}
\newcommand{\Sz}{Z.\ Szab{\'o}\,}

\newcommand{\pic}[1]{\parbox{.6cm}{\psfig{figure=#1.eps,height=.6cm}}}

\begin{document}
\title{On knot Floer homology in double branched covers}
\author[L. Roberts]{Lawrence Roberts}
\address{Department of Mathematics, Michigan State University,
East Lansing, MI 48824}
\email{lawrence@math.msu.edu}
\thanks {The author was supported in part by NSF grant DMS-0353717 (RTG)}
\maketitle

\section{Introduction}

\noindent Let $\mathbb{L}$ be a link in $A \times I$ where $A$ is an annulus. We consider $A \times I$ to be embedded in $\R^{2} \times \R$ 
respecting the obvious fibration and embedding $A$ into a round annulus in $\R^{2}$. 
We always project $\mathbb{L}$ into $\R^{2}$ (or $A$) along the $\R$-fibration. The complement of $\mathbb{L}$ in $A \times I$ is thereby
identified with the complement of $B \cup \mathbb{L}$ in $S^{3}$ where $B$ an unknot as depicted below, called the axis of $\mathbb{L}$. We assume
throughout that $\mathbb{L}$ intersects the spanning disc of $B$ in an odd number of points. For example, \\

\begin{center}
\includegraphics[scale=0.5]{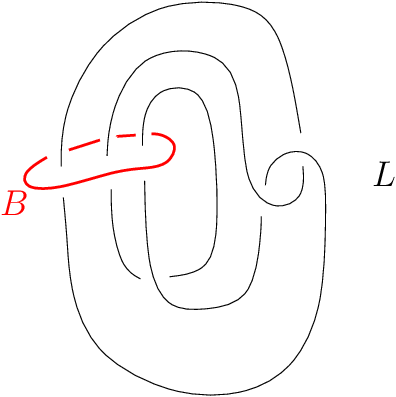}
\end{center}

\noindent Let $\Sigma(\mathbb{L})$ be the branched double cover of $S^{3}$ over $\mathbb{L}$, and let $\widetilde{B}$ be the pre-image of $B$ in $\Sigma(\mathbb{L})$. Then $\widetilde{B}$ is a null-homologous knot in $\Sigma(\mathbb{L})$ and we can try to compute 
$$
\widehat{HFK}(\Sigma(\mathbb{L}), \widetilde{B}, i) = \bigoplus_{\{\overline{\SP{s}}\, |\, \langle c_{1}(\underline{\SP{s}}), [F] \rangle\, =\, 2i\}}
\widehat{HFK}(\Sigma(\mathbb{L}), \widetilde{B}, \underline{\SP{s}}) 
$$
where $\underline{\SP{s}}$ is a relative $Spin^{c}$ structure for $\widetilde{B}$ and $[F]$ is the homology class of a pre-image of a 
spanning disc for $B$. A particularly interesting case will be when $\mathbb{L}$ is a braid. Then the pre-image of the open book of discs with binding $B$ is an open book with binding $\widetilde{B}$.\\
\ \\
\noindent To obtain a clean statement we need to adjust $\mathbb{L}$ by adding two copies of the center of $A$ which are split from the remainder of
$\mathbb{L}$. We call this new link $\mathbb{L}'$. The effect on the branched double cover is to produce $\Sigma(\mathbb{L}) \#^{2} S^{1} \times S^{2}$ containing a knot $\widetilde{B} \# B(0,0)$. We can then prove\\

\begin{prop}
Let $\mathbb{L}$ be a link in $A \times I \subset \R^{2} \times \R$ as above. Let $\mathbb{L}'$ be the adjusted version of $\mathbb{L}$. There is a spectral sequence whose $E^{2}$ term is isomorphic to the reduced Khovanov skein homology of the mirror, $\overline{\mathbb{L}}'$, in $A \times I$ with coefficients in $\mathbb{F}_{2}$ and which converges to $\oplus_{i \in \Z} \widehat{HFK}(\Sigma(\mathbb{L}) \#^{2} (S^{1} \times S^{2}), \widetilde{B} \# B(0,0), i , \mathbb{F}_{2})$. 
\end{prop}

\noindent In \cite{Doub} \Oz and \Sz constructed a spectral sequence which converged to $\widehat{HF}(Y)$ for $Y$ a double branched cover of a link in $S^{3}$. This spectral sequence featured the reduced Khovanov homology of the mirror of the link as the $E^{2}$ term. The previous proposition is a generalization of this result. \\
\ \\
\noindent In the first half of this paper, we review the skein homology, first constructed in \cite{Asae}, and examine its relationship to Khovanov homology. We then describe a spanning tree approach to computing this homology theory. This complex allows us to analyze the situation of $\mathbb{L}$ being alternating. Once this is completed we turn to building the relationship with knot Floer homology. \\
\ \\
\noindent In the second half, we derive the relationship between the two theories as the spectral sequence explained abover. We then turn to deriving some consequences of these spectral sequences. First, in \cite{Plam}, O. Plamenevskaya constructed a special element, $\widetilde{\psi}(\mathbb{L})$, of the Khovanov homology of a braid and showed that it is an invariant of the transverse isotopy class of the braid. She suggested that for certain knots, should this element survive in the spectral sequence, it would yield the contact invariant of the contact structure lifted from $S^{3}$ to the double branched cover branched over the transverse knot. This element is also a closed element in the skein homology where it defines the unique minimal filtration level. From these considerations we can prove 

\begin{prop}
Suppose there exists a $n$ such that
\ben
\item $\psi(\mathbb{L})$ is exact in the reduced Khovanov homology
\item The link surgery induced spectral sequence on $X/X_{-2g}$ collapses at $E^{2}$.
\een
then $c(\xi) = 0$.
\end{prop}

\noindent where $c(\xi)$ is the contact element for the lifted open book. The notation in this proposition is explained in section 8.\\
\ \\
\noindent Furthermore, for $\mathbb{L}$ alternating for the projection $A \times I \ra A$ much more can be said. We use the analysis of the skein homology for alternating $\mathbb{L}$ to prove the main theorem of the paper, theorem \ref{thm:alex}.

\begin{theorem}
Let $\mathbb{L}$ be a non-split alternating link in $A \times I$ intersecting the spanning disc for $B$ in an odd number of points. Then for each $k$ there is an isomorphism
$$
\widehat{HFK}(-\Si(\mathbb{L}) \#^{2} \big( S^{1} \times S^{2} \big), \widetilde{B} \# B(0,0), k) \cong \bigoplus_{i,j \in \Z} H^{i;j, 2k}(\mathbb{L})
$$
where, for each $Spin^{c}$ structure, the elements on the right side all have the same absolute $\Z/2\Z$-grading. Together these isomorphisms induce a filtered quasi-isomorphism from the $E^{2}$-page of the knot Floer homology spectral sequence to that of the skein homology spectral sequence. Thus the knot Floer spectral sequence collapses after two steps. Furthermore, for any $\SP{s} \in Spin^{c}(\Si(\mathbb{L}))$ we have that
$$
\tau(\widetilde{B}, \SP{s})  = 0
$$
where $\widetilde{B}$ is considered in $\Si(\mathbb{L})$.
\end{theorem}

\noindent In a sequel to this paper, we reprove the above theorem for $\Z$-coefficients and use it to analyze a class of fibered knots in certain three manifolds. \\
\ \\
\noindent{\bf Acknowledgements:} The author would like to thank John Baldwin and Olga Plamenevskaya for some very useful correspondence.  

\section{The Reduced Khovanov Skein Homology of \cite{Asae}}

\noindent Throughout we will assume all coefficients are in $\mathbb{F}_{2}$ and suppress the ring notation. This section gives a brief description of 
a reduced form of the theory in \cite{Asae} for categorifiying the Kauffman bracket skein module for the $I$-bundle $A \times I$ and its relationship with the reduced Khovanov homology. We adjust the account in \cite{Asae} to conform to that of Bar-Natan, \cite{BarN}. This alters the gradings from \cite{Asae} to more directly related to Khovanov's original definition. \\
\ \\
\noindent Pick an order for the $c(\mathbb{L})$ crossings in a projection, $\mathcal{P}$, of $L$
to $A$. Let $R$ be an element of $\{0,1\}^{c(\mathbb{L})}$, then associate to $R$ a collection of disjoint, simple, unoriented circles in $A$ by resolving the crossings of $\mathcal{P}$ according to: 

$$
\begin{array}{c}
\pic{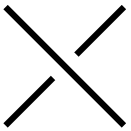}\stackrel{0}{\lra}\pic{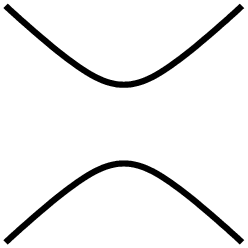} \\
\ \\
\pic{Lp}\stackrel{1}{\lra}\pic{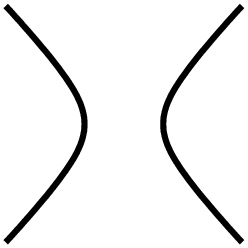}
\end{array}
$$
\ \\
\ \\
\noindent We denote the resulting diagram by $\mathcal{P}(R)$. Let $I(R)$ be  
$$
I(R) = \sum m_{i} \hspace{.15in} \mathrm{where\ } \hspace{.15in} R = \{m_{1}, \ldots, m_{n} \}  \\
$$
Finally, call an unoriented circle resulting from the resolution \textit{trivial} if it bounds a disc in $A$, and \textit{non-trivial} if it 
does not.\\ 
\ \\
\noindent An enhanced Kauffman state is then a choice of resolutions, $R$, and a choice of $\{+, -\}$ for each of the resulting circles. As usual
the enhanced states will be the generators of the chain groups. We define two bi-graded modules
$V \cong \mathbb{F} v_{+} \oplus \mathbb{F} v_{-}$ and $W \cong \mathbb{F} w_{+} \oplus \mathbb{F} w_{-}$ where
$\mathrm{deg}(v_{+}) = (1, 1)$, $\mathrm{deg}(w_{+}) = (1, 0)$ and $\mathrm{deg}(v_{-}) = -\mathrm{deg}(v_{+})$, $\mathrm{deg}(w_{-}) = -\mathrm{deg}(w_{+})$. If the resolution $R$ results in $m$ trivial circles and $l$ non-trivial circles we associate to $R$ the bi-graded module 
$$
V_{R}(\mathbb{L}) = V^{\otimes l} \otimes W^{\otimes m} \{(I(R), 0)\}
$$
We will refer to the first grading in the ordered pair as the $q$-grading and the second as the $f$-grading. \\
\ \\
\noindent The $r^{th}$ chain group, $C_{r}$ is then $\oplus_{\{R | I(R) = r\}} V_{R}(\mathbb{L})$. These will form the components of
a complex, $\mathcal{C}$, and the Khovanov skein complex will be $\mathcal{C}[-n_{-}]\{(n_{+} - 2n_{-}, 0)\}$ for some orientation on 
the link $\mathbb{L}$. The shift in $[\cdot]$ occurs in the dimension of the chain group. This last set of shifts\footnote{We follow Bar-Natan's shifting conventions} will be called the
\textit{final} shifts. We will often only be interested in relative gradings, and so will sometimes ignore the final shifts. The complex
before the final shifts will be called \textit{unshifted}. \\
\ \\
\noindent We now define the differential in the complex. As usual, we specify what happens when two circles merge in a $0 \ra 1$ resolution change, 
or what happens when a single circle divides. This suffices to specify the differential as in \cite{BarN}. The relevant maps
for merging are

$$
\begin{array}{lll}
w_{+} \otimes w_{+} \ra w_{+} & \hspace{1in} & v_{+} \otimes v_{+} \ra 0 \\
w_{+} \otimes w_{-}, w_{-} \otimes w_{+} \ra w_{-} & \hspace{1in} & v_{+} \otimes v_{-}, v_{-} \otimes v_{+} \ra w_{-} \\
w_{-} \otimes w_{-} \ra 0 & \hspace{1in} & v_{-} \otimes v_{-} \ra 0\\
\ & \ & \ \\
v_{\pm} \otimes w_{-}, w_{-} \otimes v_{\pm} \ra 0 & \ &\ \\
w_{+} \otimes v_{\pm}, v_{\pm} \otimes w_{+} \ra v_{\pm} & \ & \ \\
\end{array}
$$
\ \\
\noindent The relevant maps for dividing are 

$$
\begin{array}{lll}
w_{-} \ra w_{-} \otimes w_{-} & v_{+} \ra v_{+} \otimes w_{-} & v_{-} \ra v_{-} \otimes w_{-} \\
w_{+} \ra w_{-} \otimes w_{+} + w_{+} \otimes w_{-} & w_{+} \ra v_{+} \otimes v_{-} + v_{-} \otimes v_{+} & \ \\
\end{array}
$$
\ \\
where the rule for $w_{-}$ is determined by the topological type of the circles in the result (two trivial or two non-trivial circles).

\begin{theorem} \cite{Asae}
The tri-graded homology, $H(\mathbb{L})$, of the complex $\mathcal{C}[-n_{-}]\{(n_{+} - 2n_{-},0)\}$ with the differential defined above
is an invariant of the oriented link $\mathbb{L}$ in $A \times I$. 
\end{theorem}
\ \\
\noindent {\bf Proof:} Let $S(\mathcal{P})$ be the set of enhanced states and define for $S \in S(\mathcal{P})$

$$
\begin{array}{l}
\tau(S) =  \# \{\mathrm{positive\ trivial\ circles}\} - \# \{\mathrm{negative\ trivial\ circles}\} \\
\ \\
\Psi(S) =  \# \{\mathrm{positive\ non-trivial\ circles}\} - \# \{\mathrm{negative\ non-trivial\ circles}\}\\
\ \\
J(S) = I(S) + \tau(S) + \Psi(S) \\
\end{array}
$$
\ \\
\noindent Let $S_{ijk}(\mathcal{P})$ be the subset of $S(\mathcal{P})$ with $I(S) = i, J(S) = j,$ and $\Psi(S) = k$. Define $C^{i;jk}(\mathcal{P})$ to
be the free abelian group generated by $S_{ijk}(\mathcal{P})$.  It is shown in \cite{Asae} that the maps above define a differential on
$C^{i;jk}(\mathcal{P})$ which increases the $i$ grading by $1$. Actually, this is proved with $J'(S) = I(S) + \tau(S)$, but as the differential 
does not change $k$, the proof applies here as well. Their homology is $RII$ and $RIII$ invariant. With the shifts from a choice of orientation on the link, the theory we have outlined is also RI invariant. As with translation from Viro's notation to Bar-Natan's the shifts at the end are also necessary to pin down an invariant grading for RII, but the relative graded theory is invariant regardless. $\Diamond$\\ 
\ \\
\noindent Let $\mathcal{B}(A) \cong \{0, 1, \ldots, \} $ be the set of all link diagrams in $A$ with no crossings or trivial components. Using the 
rules

$$\pic{Lp}=\pic{L0} + tq \pic{Linfty},\quad 
L\cup \bigcirc =(q+q^{-1})L.$$
\ \\
\noindent we can associate an element of $\Z[q^{\pm 1}, t, \mathcal{B}(A)]$ to any diagram of $\mathbb{L}$, denoted $[\mathbb{L}]$. If we map the monoid $\mathcal{B}(A)$ to $\Z[q^{\pm 1}, x^{\pm 1}]$ by $1 \ra qx + q^{-1}x^{-1}$ we get a map $\phi: \Z[q^{\pm 1}, t, \mathcal{B}(A)] \ra \Z[q^{\pm 1}, t, x^{\pm 1}]$. After orienting $\mathbb{L}$, let $V(t,q,x) = t^{n_{-}}q^{n_{+} - 2n_{-}}\phi(\mathbb{L})$  which equals $\sum_{k \in \Z} q_{k, \mathbb{L}} x^{k}$ where
$$
q_{k,\mathbb{L}}= \chi_{q,t}(H^{\ast;\ast\,k}(\mathbb{L})) = \sum_{i,j}\,t^{i}q^{j} \mathrm{rk}_{\mathbb{F}}\big( H^{i;jk}(\mathbb{L})\big)
$$
 \\
\noindent The Euler characteristic for the skein homology is then $V(-1, q, x)$ and is an isotopy invariant of $\mathbb{L}$ in $A \times I$. On
the other hand $V(-1, q, 1)$ is the Jones polynomial as described by Khovanov (see also \cite{BarN}). \\
\ \\ 
\noindent There is also a reduced version of this theory.  We mark the circle in $\mathcal{P}$ that is closest to the center, at the point intersecting the spanning disc for $B$. Every diagram $\mathcal{P}(R)$ inherits this marking. Note that the marked circle in the resolved diagrams may be either trivial or non-trivial. The reduced homology is then defined to be the homology of the quotient of the above complex by the subcomplex generated by the enhanced states assigning a $-$ sign to the marked circle. The reduced chain groups are denoted $\widetilde{V}_{\mathbb{L}}(\mathcal{P})$ and the overall homology by $\widetilde{H}^{i;jk}$.

\begin{lemma}
For each $j$, there is a spectral sequence whose $E^{1}$ term is $\oplus_{i,k} H^{i;jk}(\mathbb{L})$ and which collapses at $E^{2}$  
to $\oplus_{i} H^{i,j}(\mathbb{L})$ where $H^{i,j}(\mathbb{L})$ is the usual Khovanov homology for the embedding $\mathbb{L} \ra A \times I \ra S^{3}$.
This statement also applies to the reduced theory. 
\end{lemma}
\ \\
\noindent {\bf Proof:} The entire construction has been performed so that by ignoring the distinction between trivial and non-trivial circles we obtain the Khovanov chain groups, i.e. if we use $\mathbb{L} \ra A \times I \ra \R^{2} \times I$ as an embedding of $\mathbb{L}$ in $S^{3}$ and ignore the axis. In this case we neglect the $f$-grading and treat $v_{\pm}$ and $w_{\pm}$ the same. The maps defining the differential above are almost those for the Khovanov homology, with the exception of a few terms which have been dropped. These terms are boxed below:

$$
\begin{array}{lll}
v_{+} \ra v_{+} \otimes w_{-} + \fbox{$v_{-} \otimes w_{+}$} & \hspace{1in} &
v_{+} \otimes v_{+} \ra \fbox{$w_{+}$} \\
\ \\
v_{+} \otimes w_{-}, w_{-} \otimes v_{+} \ra \fbox{$v_{-}$} & \hspace{1in} &
w_{-} \ra \fbox{$v_{-} \otimes v_{-}$} \\
\end{array}
$$
\ \\
\noindent Each of these terms preserves the $q$-grading, increases the $i$ grading by $1$, but decreases the $f$-grading by $2$. Thus, the axis can be seen as filtering the Khovanov homology, with the $E^{1}$ term of the corresponding spectral sequence being the Khovanov skein homology. Since the maps
in the spectral sequence also preserve the $-$ subcomplex, this conclusion occurs for the reduced homology as well. $\Diamond$ \\
\ \\
\begin{lemma}
Let $\mir{\mathbb{L}}$ be the mirror of $\mathbb{L}$. Then there is an isomorphism
$$
H^{i;jk}(\mathbb{L}) \cong H_{-i;-j,-k}(\mir{\mathbb{L}}) 
$$
where $H_{i;jk}$ is the corresponding cohomology group. Over a field, $\mathbb{F}$, the last group is also isomorphic to $ H^{-i;-j,-k}_{\mathbb{F}}(\mir{\mathbb{L}})$. Furthermore, the spectral sequence converging to Khovanov homology on $H^{\ast;\ast\ast}(\mathbb{L})$ is filtered chain isomorphic to that induced on the cohomology groups $H_{\ast;\ast\ast}(\overline{\mathbb{L}})$ by the 
higher differentials on $H^{\ast;\ast\ast}(\mir{\mathbb{L}})$. 
\end{lemma}
\ \\
\noindent {\bf Proof:} Each state for $\mathbb{L}$ defines
a state for $\mir{\mathbb{L}}$ by reversing the sign assignment on each circle. In addition, $0$ resolutions are now $1$ resolutions and vice-versa. Thus, $i \ra c(\mathbb{L}) - i$, $j \ra c(\mathbb{L}) - j$, and $k \ra - k$ in the unshifted theory. Examining the differential for between two states shows that the differential for $\mir{\mathbb{L}}$ is the differential for the cohomology of $\mathbb{L}$. Furthermore, after the final shifts we have $(i,j,k) \ra ( i - n_{-}, j + n_{+} - 2n_{-}, k)$ for $\mathbb{L}$ and $(c - i, c- j, - k) \ra (c - i - n_{+}, c- j + n_{-} - 2n_{+}, - k)$, where $n_{-}$ and $n_{+}$ refer to $\mathbb{L}$. This last triple equals $(-(i- n_{-}), -(j + n_{+} - 2n_{-}), - k)$. For coefficients in a field standard homological algebra implies that
$$
H_{i;jk}^{\mathbb{F}}(\mir{\mathbb{L}}) \cong H^{i;jk}_{\mathbb{F}}(\mir{\mathbb{L}})
$$
Carefully examining the terms giving rise to the spectral sequence shows that these map to the terms in the spectral sequence on the cohomology. $\Diamond$\\
\ \\
\noindent Since the $k$-grading filters the Khovanov complex, we can define for any element $\xi \in KH^{i,j}(\mathbb{L})$ a number $$T_{\mathbb{L}}(\xi) = \mathrm{min}\{k : \xi \in \mathrm{Im}\big( H_{\ast}(\oplus_{l \leq k} C^{i;jl}) \ra KH^{i,j}(\mathbb{L}) \big)\}.$$ When $\mathbb{L}$ is an unknot these numbers satisfy a relation similar to the $\tau$ invariant in knot Floer homology.   \\

\begin{lemma}
Assume $\mathbb{L}$ is an unknot and let $\mir{\mathbb{L}}$ be its mirror image. Let ${\bf u}_{\pm}$ be the generators of the Khovanov homology of the unknot in $q$-gradings $\pm 1$. Then
$$
T_{\mathbb{L}}({\bf u}_{\pm 1}) = - T_{\mir{\mathbb{L}}}({\bf u}_{\mp 1})
$$
\end{lemma}
\ \\
\noindent {\bf Proof:}  Let $\mathcal{F}_{j;s} = \oplus_{i; k \leq s} C^{i;jk}(\mathbb{L})$ and let $C_{j} = \oplus_{i,k} C^{i;jk}(\mathbb{L})$. Since the differential preserves the $q$-grading, $j$, there is a long exact sequence:
$$
0 \lra \mathcal{F}_{j;s} \stackrel{I_{s}}{\lra} C_{j} \stackrel{P_{s}}{\lra} Q_{j;s} \lra 0
$$
where $Q_{j;s}$ is the quotient complex, $C_{j}/ \mathcal{F}_{j;s}$. Now $\oplus_{j} H_{\ast}(C_{j}) =\Z{\bf u}_{+} \oplus \Z{\bf u}_{-}$, and $T_{\mathbb{L}}$ measures the first $s$ for which the map in the long exact sequence on homology will include ${\bf u}_{\pm}$ in the image of $I_{s\ast}$ relative to the $q$-grading. \\
\ \\
\noindent There is a duality isomorphism $D: H^{i;j}(U) \ra H_{-i;-j}(\mir{U})$, $D({\bf u}_{\pm}) = {\bf u}_{\mp}$, on the Khovanov homologies which is induced by the symmetric pairing $a_{+} \otimes a_{-} \stackrel{m}{\lra} a_{-} \stackrel{\epsilon}{\lra} 1$ where $\epsilon: A \ra \Z$ is the counit for the Frobenius algebra underlying Khovanov homology. In particular, $a_{+} \ra \langle a_{+}, \cdot \rangle = a_{-}^{\ast}$. This can be extended to $V$ as well, and corresponds to changing the markers on each of the circles in an enhanced state. It thus induces a map on the skein homology spectral sequences. Checking the effect on the differential establishes the following commutative square:

$$
\begindc{\commdiag}[6]
	\obj(0,0)[o1]{$0$}
	\obj(0,10)[o2]{$0$}
	\obj(10,0)[I1]{$Q^{\ast}_{-1;-s-1}(\mir{U})$}
	\obj(25,0)[I2]{$C^{\ast}_{-1}(\mir{U})$}
	\obj(10,10)[I3]{$\mathcal{F}_{+1,s}(U)$}
	\obj(25,10)[I4]{$C_{1}(U)$}
	\obj(40,0)[I5]{$\mathcal{F}^{\ast}_{-1,-s-1}(\mir{U})$}
	\obj(40,10)[I6]{$Q_{+1;s}(U)$}
	\obj(50,0)[o3]{$0$}
	\obj(50,10)[o4]{$0$}
	\mor{o1}{I1}{$\,$}
	\mor{o2}{I3}{$\,$}
	\mor{I1}{I2}{$P_{-s-1}^{\ast}$}
	\mor{I3}{I4}{$I_{s}$}
	\mor{I4}{I2}{$D$}
	\mor{I3}{I1}{$D$} 
	\mor{I2}{I5}{$I_{-s-1}^{\ast}$}
	\mor{I4}{I6}{$P_{s}$}
	\mor{I6}{I5}{$D$}
	\mor{I6}{o4}{$\,$}
	\mor{I5}{o3}{$\,$}
\enddc 
$$
\ \\
\noindent If ${\bf u}_{+}$ is in the image of $I^{s}$, then ${\bf u}_{-}^{\ast}$ is in the image of $P^{\ast}_{-s-1}$. In paticular, $I^{\ast}_{-s-1}({\bf u}_{-}^{\ast}) = 0$. But then there is no element in $\mathcal{F}_{-1,-s-1}$ which maps to ${\bf u}_{-}$ and so $-s - 1 < T_{\mir{\mathbb{L}}}({\bf u}_{-})$. Thus if $s = T_{\mathbb{L}}({\bf u}_{+})$ then $-T_{\mathbb{L}}({\bf u}_{+}) \leq T_{\mir{\mathbb{L}}}({\bf u}_{-})$. If ${\bf u}_{+}$ is not in the image of $I^{s}$ -- i.e. $s < T_{\mathbb{L}}({\bf u}_{+})$ -- then $I^{\ast}_{-s - 1}({\bf u}_{-}^{\ast}) \neq 0$. Choose some element on which this image pairs non-trivially and is uniformly in $q$-grading $-1$. This element must then map in homology to $a_{-}$ and $-s - 1 \geq T_{\mir{\mathbb{L}}}({\bf u}_{-})$. In particular, if $s = T_{\mathbb{L}}({\bf u}_{+}) - 1$ then $- T_{\mathbb{L}}({\bf u}_{+}) \geq T_{\mir{\mathbb{L}}}({\bf u}_{-})$. This proves the result.  
$\Diamond$\\
\ \\
\noindent Let $\mathbb{L}_{1}$ and $\mathbb{L}_{2}$ be two links in $A \times I$. Let $\mathbb{L} = \mathbb{L}_{1} | \mathbb{L}_{2}$ be the link in 
$A \times I$ where $A = \{z : 1 \leq |z| \leq 3\}$ and $\mathbb{L}_{1}$ lies in $\{z : 1 \leq |z| \leq 2\} \times I$ while $\mathbb{L}_{2}$ lies in 
$\{z: 2\leq |z| \leq 3\} \times I$. Then we can prove\\

\begin{lemma}
With coefficients in a field, $\mathbb{F}$, there is an isomorphism
$$
H^{i;jk}(\mathbb{L}) \cong \bigoplus_{\begin{array}{c} i_{1} + i_{2} = i,\ j_{1} + j_{2} = j \\ k_{1} + k_{2} = k \end{array}} H^{i_{1};j_{1}k_{1}}(\mathbb{L}_{1}) \otimes H^{i_{2};j_{2}k_{2}}(\mathbb{L}_{2})
$$
In fact, this is an isomorphism of spectral sequences, so that if $\xi_{1} \in KH^{i_{1};j_{1}}(\mathbb{L}_{1})$ and $\xi_{2} \in KH^{i_{1};j_{1}}(\mathbb{L}_{1})$ then 
$$
T_{\mathbb{L}}(\xi_{1} \otimes \xi_{2}) = T_{\mathbb{L}_{1}}(\xi_{1}) + T_{\mathbb{L}_{2}}(\xi_{2})
$$
\end{lemma}
\ \\
\noindent Finally, if there is a non-trivial component, $L_{1}$, split from the rest of $\mathbb{L}$, it can be made to lie in the diagram without
crossing any other strand of $\mathbb{L}$. $L_{1}$ survives unchanged in every resolution; thus, marking it induces a marking on a non-trivial circle for every resolution  When the number of intersections of $\mathbb{L}$ with the spanning disc for $B$ is odd, then the reduced skein homology of this configuration has the form $\widetilde{H}(\mathbb{L} - L_{1}) \otimes V$. This choice shifts the complex by $\{(1,1)\}$, so we will always shift at the end to compensate. Thus, the final shift will be $[-n_{-}]\{(n_{+} - 2n_{-} - 1, -1)\}$ for this marking convention. This component is special and doesn't participate in the calculations. Thus for the mirror, we keep it labelled $+$ and derive the same duality for the mirror image as previously. 
 
\section{Spanning Tree Complex}

\noindent As with Khovanov homology, the skein homology for links in $A \times I$ with connected projections admits another presentation in terms
of the spanning trees for the knot diagram. We follow \cite{Wehr} in establishing this result, but see also \cite{Cham}.\\
\ \\
\noindent Start with the projection of $\mathbb{L}$ in the plane and follow Wehrli's algorithm. First, number the crossings. Now proceed
to resolve the first crossing if both resolutions produce connected diagrams. The homology can be shown to be a mapping cone on these
resolutions. If resolving the first crossing disconnects the diagram for one or other resolution, proceed to the next crossing until you come to one where both resolved diagrams are connected or you run out of crossings. Now iterate this procedure on the resolved diagrams. The result is
a tree of diagrams, the resolution tree, whose leaves correspond to unknots that are reducible to the standard unknot using only the first Reidemeister move. There is a unique way to smooth the remaining crossings to get an unknot in the plane. Likewise to each complete smoothing which produces an unknot in the plane there is a unique leaf which smooths to it (due to the enumeration of the crossings). Let $K_{1}(\mathbb{L})$ be the complete smoothings with only one component and for each $S \in K_{1}(D)$ let $D_{S}$ be the twisted unknot corresponding to $S$.\\
\ \\
\noindent We can prove that the \textit{unshifted} skein complexes behave in the following way due to an RI move:\\

\begin{figure}
\begin{center}
\includegraphics[scale=0.5]{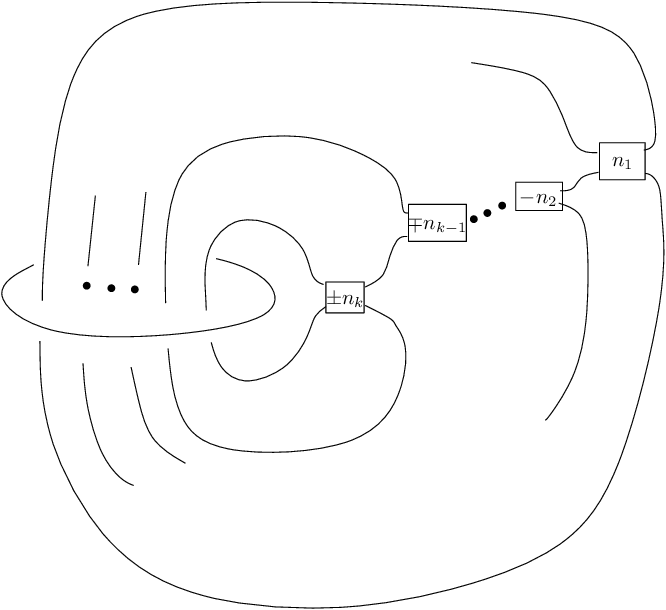}
\end{center}
\caption{The special class of unknots which act as the base cases for the spanning tree complex.} \label{fig:special}
\end{figure}

\begin{center}
\includegraphics[scale=0.75]{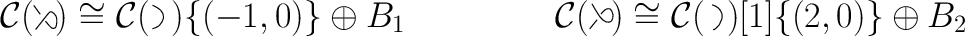}
\end{center}

\noindent where $B_{1}$ and $B_{2}$ are contractible. It does not matter whether the RI move involves trivial or non-trivial
components in the complete smoothings. The shifts result from the invariance of the theory after the final shifts are performed. 
For the first RI move above, the left side would need to be shifted $[0]\{(1,0)\}$ further than the right side, due to the 
extra positive crossing. Thus the right side should be shifted $\{(-1,0)\}$ to correspond to the left side in the unshifted complex. \\  
\ \\
\noindent With this observation we may proceed analogous to \cite{Wehr} to obtain the following proposition:

\begin{lemma} \cite{Wehr}
Let $\mathbb{L} \subset A \times I$ have a connected diagram in $A$. Then there is a decomposition $\overline{C} \cong A \oplus B$ where
$\overline{C}$ is the unshifted version of the skein complex, $B$ is contractible, and $A$ is given by
$$
A = \bigoplus_{S \in K_{1}(\mathbb{L})} H^{\ast;\ast\ast}(D_{S})[-w(D_{S})]\{(-2w(D_{S}), 0)\}[r(S)]\{(r(S), 0)\}
$$
where $w(D_{s})$ is the writhe of $D_{s}$, and $r(D,S)$ is the number of $1$ smoothings necessary in resolving $\mathbb{L}$ to get $S$.
\end{lemma}
\ \\
\noindent Of course, the unshifted homology $H^{\ast;\ast\ast}(D_{S})$ does not change while using Reidemeister I moves, but the difference between considering $D_{S}$ in $A$ and in the plane is precisely to disallow RI moves which would need to cross $B$. Thus $D_{S}$ can be simplified to $D'_{S}$, a twisted unknot, where all the twisting ultimately must link with $B$. This implies that $D'_{S}$ is isotopic to a knot of the special form in Figure \ref{fig:special}, where each $n_{i}$ records the number of half twists. The homologies of these knots form the building blocks of the spanning tree complex above.\\

\noindent Since we assumed that the number of intersections with the spanning disk is odd, the diagram $D_{S}$ must be non-trivial, and thus link $B$.
Suppose a connected diagram, fully reduced in $A \times I$ does not look as above. Choose a region ``under'' $B$ and follow it around clockwise and counter-clockwise. Suppose in both directions we come to crossing regions which traverse the region, as we must since the diagram is connected. Suppose
that these twist regions are distinct. Since the diagram is fully RI reducible, only the strands of one of the twist regions can link $B$. The other
must then reduce using RI moves, since the diagram is essentially planar. This forces a knot isotopic to one of the unknots seen above.  \\
\ \\
\noindent We can say a little more concerning the unknots in \ref{fig:special}. In particular, we compute the numbers, $T_{\mathbb{L}}$, for these unknots. Let $\#(T_{\pm})$ be the number of left/right-handed twist regions in Figure \ref{fig:special}. Since these are unknots, their Khovanov homologies are composed of $\mathbb{F}{\bf u}_{+}$ in homological and $q$-grading $(0,1)$ and $\mathbb{F}{\bf u}_{-}$ in $(0,-1)$. Finally, assume that $\mathbb{L}$ links the axis an odd number of times. We can then prove:\\

\begin{prop}
For the special unknotted branch loci in Figure \ref{fig:special}, let $\,T(\mathbb{L})$ denote $\#(T_{-}) - \#(T_{+})$. Then $$T_{\mathbb{L}}({\bf u}_{\pm 1}) = T(\mathbb{L}) \pm 1.$$ For the alternating unknots in this family, $T_{\mathbb{L}}({\bf u}_{\pm 1}) = \pm 1$.
\end{prop}
\ \\
\begin{figure}
\begin{center}
\includegraphics[scale=0.5]{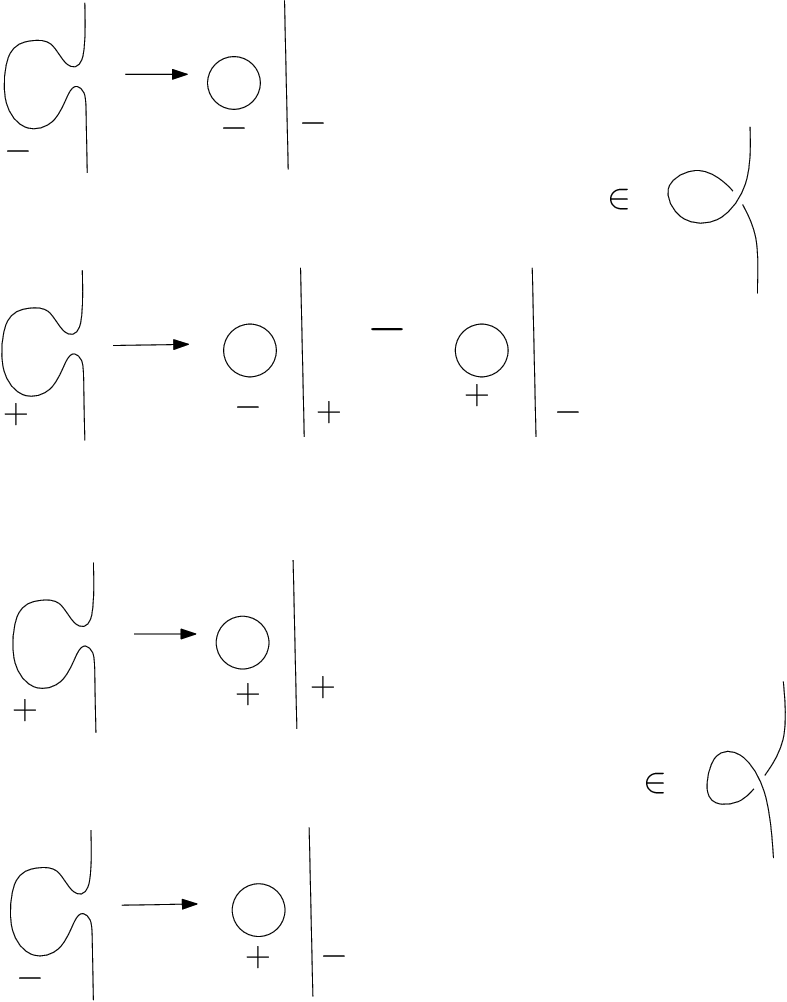}
\end{center}
\caption{Rules for transfering generators when an RI move is applied. The particular twist is represented on the far right. Note that these maps are chain maps inducing isomorphisms on the Khovanov homologies, \cite{Khov}} \label{fig:RIred}
\end{figure}

\noindent {\bf Proof:}  These unknots are isotopic to the standard planar unknot using only $RI$-moves. M. Jacobsson provides rules for mapping closed elements in the Khovanov cube of a link to those of the link with a single $RI$ move, which in our notation are as in Figure \ref{fig:RIred}. We can use these moves to try to compute $T_{\mathbb{L}}({\bf u}_{\pm 1})$. As a first step, we exhibit a specific generator which will produce ${\bf u}_{\pm 1}$ in homology. The maximal value of $k$ needed to obtain this generator in $\oplus_{l \leq k} C^{i;jk}$ will then be an upper bound on $T_{\mathbb{L}}({\bf u}_{\pm 1})$. \\
 \ \\
 
\begin{figure}
\begin{center}
\includegraphics[scale=0.5]{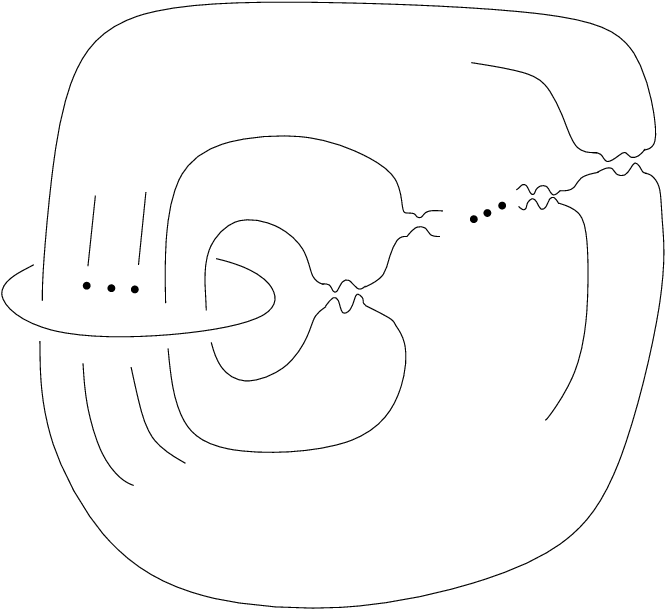}
\end{center}
\caption{ } \label{fig:Start}
\end{figure} 
 
 \noindent Consider an unknot as in Figure \ref{fig:Start} formed by resolving all the crossings of $\mathbb{L}$ horizontally.  
We begin by examining the effect of replacing the outermost resolution with a crossing, according to the Jacobbson rules. According to these rules,  for left handed twist regions, a $-$ on the outer circle will propagate to each new circle as we change the resolution at each crossing. Meanwhile
a $+$ maps to a linear combination of the two generators formed by a $+$ on one of the new circles and a $-$ on the other, or vice-versa.
For right handed twist regions the argument is slightly different. Namely, $+$ markers are placed on the inner circles regardless of the marker on the outer circle.  This is a chain map since the disconnected diagram occurs as the $1$-resolution for the new crossing.   \\
 \ \\
\noindent When resolving the outermost crossing there are four cases to consider:
\ben
\item The crossing occurs in a left handed twist region, and the original circle is labelled with a $-$. Then each new circle will be labelled with a minus until we come to a right handed twist region. All the non-trivial circles formed by the resolution changes will be labelled with a minus and thus
we have a contribution of $-\#(T_{+}) - 1$ to the $k$-grading, as this is the number of non-trivial circles in this group. The right handed twist region which may follow will be of the type (2), to which we turn now. 
\item The crossing is in a right handed twist region, and the circle is labelled with a $-$. Then all the new circles will be labelled with $+$'s until the next left handed region, which is of type (3). There are $\#(T_{-})$ nontrivial circles which receive a $+$ marking in this way. If the $-$ is the outermost in the whole diagram, we have a contribution of $\#(T_{-}) - 1$. If it comes from a left handed region preceding, say from case (1), then we have added $\#(T_{-})$ to the amount already there, which we assume satisfies the proposition. Thus, we will still satisfy the proposition after this right handed region, especially if it is the last. For example, if we pass from type (1) to type (2), we can associate each marker with the twist region on its right, including here the unbounded complement of the diagram as a left handed region. This gives a total of $\#(T_{-}) - \#(T_{+}) - 1$.
\item The crossing is in a left handed region, and the circle is labelled with a $+$. This is the case where we need the linear combination of generators. The crucial observation is that the linear combination, which grows at each crossing in the twist region, ultimatley involves generators
with at most one $+$ marker abutting the crossings in the consecutive left handed twist regions. If we place the $+$ on a circle and the $-$ on the new
circle, then all additional circles until we change handedness will have a $-$. However, if we place the $+$ on the new circle, and a $-$ on the old
circle, we will generate a string of $-$'s to the right of the $+$. If the $+$ marker is on a non-trivial circle, the total number of non-trival plus circles does not change, whereas there will be $\#(T_{+})$ minus markers introduced. If it is on a trivial circle, then there are $\#(T_{+}) + 1$
minus markers introduced and these generators are in a {\it smaller} $k$-grading. Whatever marker winds up on the non-trivial circle at the junction with the next right handed twist region is immaterial as both type (2) and (4) will propagate the same number of additional $+$ markers.   
\item If the crossing is in a right handed twist region, and the circle is labelled with a plus we obtain plus markers on all the new circles until the
next left handed twist region. If the circle is the outermost, this is a contribution of $\#(T_{-}) + 1$. If it follow another region, we have added
$\#(T_{-})$ to the running total, and the proposition is still satisfied.
\een
Note that type (1) can only occur in the very outermost twist region, since the right handed twist regions always pass a $+$ marker to the next left handed region. It is type (3) which truly determines the outcome. Checking the numerics shows that the maximal $k$-grading for the generators in the linear combination so produced is $\#(T_{-}) - \#(T_{+}) \pm 1$ where the sign of $\pm 1$ is determined by whether we start with a $+$ or $-$ marker on the unknot in Figure \ref{fig:Start}. Since these are generators of the Khovanov homology for the unknot, and the maps in Figure \ref{fig:RIred}
are the chain maps used to show the RI-invariance of Khovanov homology, we have exhibited an element of the chain complex of the skein homology which survives the spectral sequence and will represent one of ${\bf u}_{\pm}$ depending on the original marker.\\
\ \\
\noindent Altogether, this shows that $T_{\mathbb{L}}({\bf u}_{\pm 1}) \leq \#(T_{-}) - \#(T_{+}) \pm 1$. However, the argument also applies to $\overline{\mathbb{L}}$ and we know that $T_{\mir{\mathbb{L}}}({\bf u}_{\mp 1}) = - T_{\mathbb{L}}({\bf u}_{\pm 1})$. In the mirror image there are $\#(T_{+})$ left handed regions and $\#(T_{-})$ right handed regions. Hence, $T_{\mir{\mathbb{L}}}({\bf u}_{\mp 1}) \leq \#(T_{+}) - \#(T_{-}) \mp 1$. Replacing  the left side with $- T_{\mathbb{L}}({\bf u}_{\pm 1})$ gives $T_{\mathbb{L}}({\bf u}_{\pm 1}) \geq \#(T_{-}) - \#(T_{+}) \pm 1$, and the result follows The final statement is simply a reflection of the even number of twist regions, alternating between handedness, when there are an odd number of strands. $\Diamond$

\section{Results for the skein homology of alternating links}

\noindent The goal of this section is to use the spanning tree presentation of the skein homology to prove the following theorem

\begin{theorem}
Let $\mathbb{L}$ be an alternating link in $A \times I$ intersecting the spanning disc for $B$ in an odd number of points. Then the Khovanov
skein homology $H^{i;jk}(\mathbb{L})$ is trivial unless $k - j + 2i = \sigma(\mathbb{L})$. Thus the homology is determined by the Euler generating polynomial $V(t, q, x) = t^{n_{-}}q^{n_{+} - 2n_{-}}\phi([\mathbb{L}])$, defined in section 1, and the signature of the oriented link $\sigma(\mathbb{L})$, thought of as embedded in $S^{3}$.   
\end{theorem}

\noindent We will follow \cite{Wehr} in calculating the Khovanov-type homology of an alternating configuration. Both provide simplified proofs of E. S. Lee's result concerning alternating links, \cite{ELee}, which describes the result of computing the spectral sequence for the axis filtration: the homology will be supported on the lines $j - 2i = - \sigma(L) \pm 1$. It is towards a variation of this result that we now aim. Note, however, that our result is not just about supports. We return to this at the end of the proof. \\
\ \\ 
\noindent Assume that $\mathbb{L}$ admits an alternating projection to $A$ which is connected as a subset of $A$. We maintain the assumption that $\mathbb{L}$ intersects the spanning disc for $B$ in an odd number of points; however, we will relax this when it is to our advantage. We will bi-color the \textit{plane} according to the following convention:

\begin{center}
\includegraphics[scale=0.4]{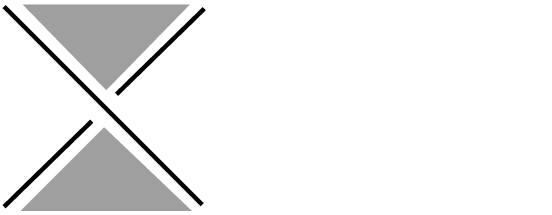}
\end{center}

\noindent For any $\mathbb{L}$, regardless of the parity of intersecting the spanning disk, we define $M(\mathbb{L})$ to be the number $N_{W} - N_{B}$
where $N_{W}$ is the number of white regions intersecting the projection of $B$ and $N_{B}$ is the number of black regions. When $\mathbb{L}$ 
intersects the spanning disc in an odd number of points, $M(\mathbb{L}) = 0$; for an even number of points $M(\mathbb{L}) = \pm 1$. This number
does not change under Reidemeister moves applied to $\mathbb{L}$, nor does it change when crossings of $\mathbb{L}$ are resolved. Furthermore, all
the projections, $D_{S}$, in the spanning tree complex will be alternating. We start with a lemma concering these unknots \\

\begin{lemma}
For each alternating twisted unknot in Figure \ref{fig:special} the homology $H^{i;jk}(\mathbb{L})$ satisfies
$k - j + 2i = M(\mathbb{L})$. 
\end{lemma}

\noindent We will show that diagrams of the special form above have the property that $H^{i;jk}$ satisfies $k - j + 2i = M(D_{S})$, and that the last number is determined by the type of crossing on the outermost boundary.  We start with the following cases:

\ben
\item $\mathbb{L}$ as a single non-trivial unknot has this property. Its homology is $0$ unless $(i;j,k) = \pm (0;1,1)$, and those
have homology $\mathbb{F}$. But then $k - j + 2i = 0 = M$ since there is one black and one white region.
\item If $D_{S}$ has the property that $k - j + 2i = C$ so does $D_{S} \cup N$ where $N$ is a disjoint non-trivial circle. 
\item The closures of $\sigma_{1} \in B_{2}$ and $\sigma_{1}^{-1} \in B_{2}$ have the property that $k - j + 2i = M(D)$. This requires a computation.
For $\sigma_{1}^{-1}$ the shifted complex has homology 
$$
H^{i;jk} \cong \left\{ \begin{array}{cl} \mathbb{F}_{-1} & (j,k) = (-3,0) \\
																				\mathbb{F}_{0} & (j,k) = (-3,-2), (-1, 0), (1, 2) \\ \end{array} \right. 
$$
where the subscript denotes $i$, and each element has $k - j + 2i = + 1$. Furthermore, $N_{W} = 2$ and $N_{B} = 1$, so $M(D) = 1$. For the closure
of $\sigma_{1}$ we obtain:
$$
H^{i;jk} \cong \left\{ \begin{array}{cl} \mathbb{F}_{1} & (j,k) = (3,0) \\
																				\mathbb{F}_{0} & (j,k) = (-1,-2), (1, 0), (3, 2) \\ \end{array} \right. 
$$
and $k - j + 2i = -1 = M(D)$.
\een

\begin{figure}
\begin{center}
\includegraphics[scale=1.0]{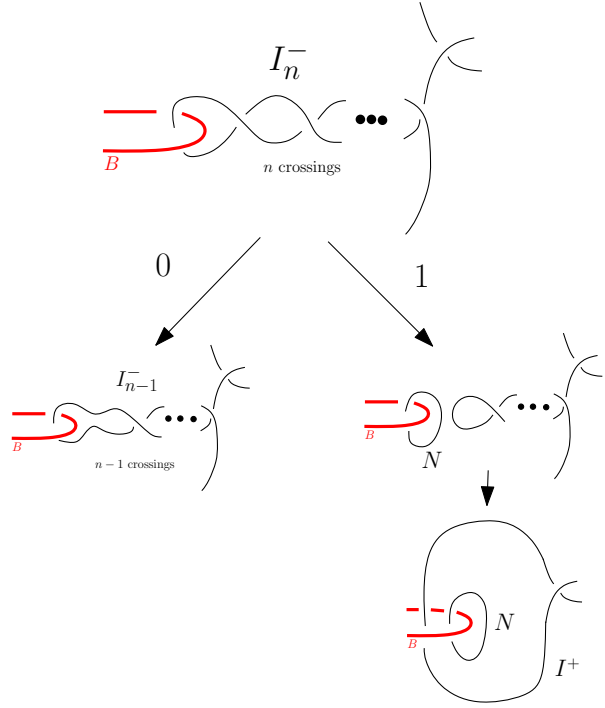}
\end{center}
\caption{A depiction of $I^{-}_{n}$ for $n > 1$, and the corresponding $I^{+}$, as it occurs in the resolution tree for the innermost crossing}. \label{fig:Iminus}
\end{figure}

\begin{figure}
\begin{center}
\includegraphics[scale=0.6]{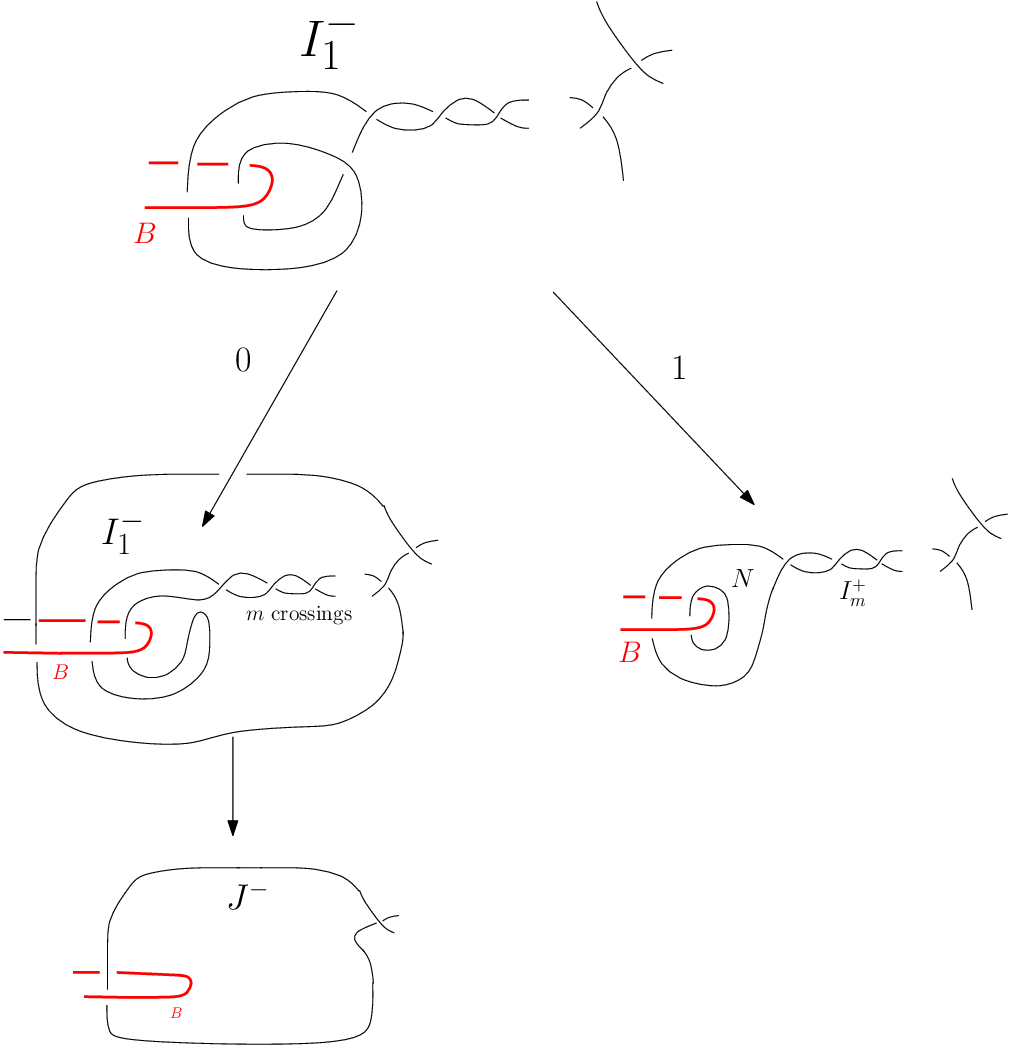}
\end{center}
\caption{A depiction of $I^{-}_{1}$, $I^{+}$ and $J^{-}$, as they occur in the resolution tree for the innermost crossing}\label{fig:Jminus}
\end{figure}

\noindent The nontrivial unknot, and the closures of $\sigma_{1}$ and $\sigma_{1}^{-1}$, are the base cases for our induction. We now assume that we have a twisted unknot of the special type above, which has $\sigma(\mathbb{L}) = 0$. See Figure \ref{fig:Iminus} and Figure \ref{fig:Jminus} to clarify the notation. We start by assuming that near the inner point where $B$ crosses the plane the twisting is right-handed. Assume that there are $n > 1$ negative crossings, and call this knot $I_{n}^{-}$. If we $0$ resolve the innermost crossing we obtain $I^{-}_{n-1}$, \\
\ \\
while if we $1$ resolve the crossing we obtain $N \cup I^{+}$. Let $[s]\{(t,0)\}$ be the contribution to the final shift of $I_{n}^{-}$ arising
from the crossings not involoved in this twist region. There is then a long exact sequence
$$
\begin{array}{c}
\lra H^{\ast}(I^{+} \cup N)[-s]\{(-t,0)\}[n-1]\{(2n-2,0)\}[1]\{(1,0)\} \lra \ \\ \\ H^{\ast}(I_{n}^{-})[-s + n]\{(-t + 2n,0)\} \lra 
H^{\ast}(I_{n-1}^{-})[-s + n-1]\{(-t + 2n-2,0)\} \lra\\
\end{array}
$$
where the sequence arises from $0 \ra 1$ resolution maps in the unshifted complexes. The additional shifts for $H^{\ast}(I^{+}\cup N)$ come from its arising in the $1$ resolution and from the additional negative crossings introduced from the twists needed before the resolution change. Those for $H^{\ast}(I_{n-1}^{-})$ come from the negative crossings remaining in the resolved diagram. The maps two internal maps are degree preserving. If
$k - j + 2i = C$ for the shifted complex for $I_{n}^{-}$ then $k- j + 2i = C + t - 2n + 2 -2s + 2n - 2 = C + t - 2s$. Reversing the shift for $H^{\ast}(I_{n}^{-})$ we obtain that elements mapping to the homology for $I_{n-1}^{-}$ satisfy $k - j + 2i = C + t - 2s - 2n + t + 2n + 2s = C$.
By assumption, $C = M(I_{n-1}^{-})$ and $M(I_{n}^{-}) = M(I_{n-1}^{-})$ since there has been no change in the black/white region count. On the other hand, $M(I^{+}) = M(I_{n}^{-}) - 1$ since we have lost the interior region, necessarily white by our crossing assumption. Since the addition of $N$
does not change $k - j + 2i$, if $k- j + 2i = M(I^{+})$ for $I^{+}$, we see that the terms in the unshifted complex for $I^{+}$ satisfy
$k - j + 2i = M(I^{+}) + t - 1 - 2n + 2 + 2(-s + 1 + n - 1)$ $= M(I^{+}) + t - 2s + 1$. Applying the shifts to get the shifted complex for $I_{n}^{-}$
prodeuces elements with $k - j + 2i = M(I^{+}) + 1 = M(I_{n}^{-})$. Every element in the image of $H^{\ast}(I^{+} \cup N)$ will also have the desired property. Thus by induction, the property will be true also for $I_{n}^{-}$. \\
\ \\
\noindent This leaves the case where $n=1$. The $1$ resolution occurs in the same way and we may draw the same conclusion. However, for the $0$
resolution a large collapse can occur. If $I^{+}$ has $m \geq 1$ positive crossings in the next region, the $0$ resolution allows us to untwist all 
of these until we get to $J^{-}$. The complex for $J^{-}$ is thus shifted by $\{(-m, 0)\}$ when injected into that for $I^{-}_{1}$. This implies that
$k- j + 2i$ increases by $m$ in the unshifted complexes. In the shifted complexes, $I^{-}_{1}$ is shifted $[-1]\{(m - 2, 0)\}$ more than $J^{-}$. That
shift reduces $k - j + 2i$ by $2-m - 2 = - m$. Thus after the final shift there is a difference of $0$. But note that the resolution eliminates both
a black and a white region and thus leaves $M(J^{-}) = M(I^{-}_{1})$. All told, if $k - j + 2i = M$ holds for the knots with fewer crossings and the innermost crossing is negative then it also holds for $I^{-}_{n}$.\\
\ \\
\noindent A similar argument can be deployed for the case where the innermost crossing is positive. Alternately we can appeal to the symmetry under 
reflection to switch the two cases. Since this switches the black and white regions, it also multiplies $M$ by $-1$. $\Diamond$\\
\ \\
\noindent Thus for every unknot in our collection we have $k - j + 2i = M(K)$ for every generator in the homology. In particular, $(j,k)$ determines $i$.
Note that this conclusion remains valid if we add a single marked non-trivial circle. It also remains true if we shift by $[-w]\{-2w\}$. As with the original proofs of the alternating links property, the value of $r(S)$ is the same for every complete smoothing in $K_{1}(\mathbb{L})$, depending
only on the number of black regions and the crossings joining them.  So all the generators for the spanning tree model of the unshifted homology satisfy $k - j + 2i = r(S)$ after the $[r(S)]\{(r(S), 0)\}$ shifts and the odd number of intersections. 
The final shift of the diagram for $\mathbb{L}$ is $[-n_{-}]\{(n_{+} - 2n_{-}, 0)\}$ and produces generators satisfying $k - j + 2i = r(S) - n_{+}$. 
From \cite{ELee}, we have that $r(S) - n_{+} = \sigma(\mathbb{L})$. Thus, after the final shifting, every generator in
the spanning tree complex satisfies $k - j + 2i = \sigma(\mathbb{L})$ and since the differential in the spanning tree complex also preserves $(j,k)$ and increases $i$, this is also the homology. For those generators which survive the spectral sequence to the 
Khovanov homology, we also have that $j - 2i = -\sigma(\mathbb{L}) \pm 1$. Thus, for these generators, $k = \pm 1$. \\
\ \\
\noindent {\bf A comment about supports:} Wehrli's argument produces an unshifted chain complex which has the same chain groups for $l + r = i$ and $2l + r \pm 1 = j$ where $r= r(S)$ is constant. Thus $j - 2i = - r \pm 1$ which when shifted yields $j - 2i = -\sigma(\mathbb{L}) \pm 1$. For a given $q$-grading, $j$, there are two $i$ gradings differing by $1$. Thus there can still be non-zero terms in the differential, which may result in
torsion or vanishing homology groups, and thus the homology is at most supported on these lines. In our case, these groups are distinguished by their
$k$-value, which is also preserved by the differential. The issue of torsion returns in the spectral sequence, but it is known that at most $2^{r}$-torsion occurs for alternating knots, \cite{Shum}, and so working over $\mathbb{F}_{2}$ will correct it.         
 
\section{Resolutions in knot Floer homology}

\noindent We now leave the Khovanov skein homology to establish some results linking it to knot Floer homology. The two will intertwine in later sections. \\
\ \\
\noindent Assume that $\mathbb{L}$ intersects a spanning disc for $B$ generically in an odd number of points. Let $R$ be a complete resolution of the crossings in $\mathcal{P}$. Of the closed curves in $\mathcal{P}(R)$, some number, $m$, are geometrically split from the axis, $B$. The remainder, $l$, form an unlink each of whose components link the axis one time. For such a link of unknots, the double branched cover is easily computed to be $\#^{l + m - 1} S^{1} \times S^{2}$. Moreover, $\widetilde{B}(R)$ is still a knot since each unknot which is split from $B$ intersects a disc generically an even number of times. This knot is $\#^{\frac{l -1}{2}} B(0,0)$ $\subset \#^{l + m - 1} S^{1} \times S^{2}$ (and the unknot in $\#^{m} S^{1} \times S^{2}$ if $l=1$) where $B(0,0) \subset S^{1} \times S^{2} \# S^{1} \times S^{2}$ is the knot obtained by performing $0$ surgery on any two of the three components of the Borromean rings. Hence, 
 $$
 \widehat{HFK}(\widetilde{B}) \cong V^{\otimes(l-1)} \otimes W^{\otimes m}
 $$
 \ \\
\noindent where $V \cong \Z_{(\frac{1}{2}, \frac{1}{2})} \oplus \Z_{(-\frac{1}{2}, - \frac{1}{2})}$ and $W \cong \Z_{(\frac{1}{2}, 0)} \oplus \Z_{(-\frac{1}{2}, 0)}$. Here the first term in the subscript is the rational grading, whereas the second term is the filtration. Since $l - 1$ is even, the filtration levels are in fact integers. Furthermore, the homology is entirely supported in the trivial $Spin^{c}$ structure.  Note that there are no higher differentials in the knot Floer spectral sequence. All we have done is compartmentalize in a new manner the Heegaard-Floer homology of $\#^{l + m - 1} S^{1} \times S^{2}$. It is the latter which is associated to $R$ in \cite{Doub}. \\
\ \\
\noindent We wish to define an isomorphism 

$$
\Phi_{B}(R) : \widetilde{V}_{\mathbb{L}}(R) \stackrel{\cong}{\lra} \widehat{HFK}(\Sigma(\mathcal{P}(R)), \widetilde{B})
$$
\ \\
\noindent for any complete resolution $R$. However, a slight mismatch arises: the knot Floer homology of the binding
implicitly corresponds to marking a non-trivial circle. This cannot always be arranged in the skein homology theory. We rely upon a 
trick to resolve this problem: we introduce two non-trivial circles into $\mathbb{L}$ which link $B$ once and otherwise do not
interact with the diagram. These should be considered innermost circles. We always mark the innermost one (we need two to keep the binding
connected) and since this circle does not include any crossings, it will be the marked circle throughout. \\
\ \\
\noindent The effect on the double cover of changing $\mathbb{L}$ to $\mathbb{L}'$ is to replace $\Sigma(\mathbb{L})$ with $\Sigma(\mathbb{L}) \#^{2} S^{1} \times S^{2}$ and to replace $\widetilde{B}$ with $\widetilde{B} \# B(0,0)$. We see this by shrinking the two new components to nearby meridians of $B$ and then examining the double cover of a small ball which includes them and an arc on $B$. The effect of these connect sums on the Heegaard-Floer homology is well understood. In particular, since $B(0,0)$ induces an entirely collapsed spectral sequence for the Heegaard-Floer homology, we will be able to read off any information about the knot Floer homology of $\widetilde{B}$ from that of $\widetilde{B} \# B(0,0)$. It would be nice to avoid the introduction of new components, but doing so only adds a slight increase in complexity. \\
\ \\
\noindent With this alteration, order the circles in $\mathcal{P}(R)$ by the marked circle first, then all the non-trivial circles, then all the trivial circles. An element of $\widetilde{V}(\mathcal{P}(R), B)$ is encoded as $+ \otimes v^{1}_{\pm} \otimes \cdots \otimes w^{n}_{\pm}$ and is mapped to $\gamma_{i_{1}}\cdots \gamma_{i_{k}} \cdot \Theta^{+}$ where $\{ i_{1}, \ldots, i_{k} \}$ are the indices for the minus signs on non-marked circles, $\gamma_{j}$ is the first homology class dual to the $j^{th}$ sphere, and $\Theta^{+}$ is the highest degree generator of $\widehat{HF}(\Sigma(\mathcal{P}(R)))$. In particular, a representative for $\gamma_{j}$ in $\mathbb{F}_{2}$-homology can be found by lifting an arc between the marked circle and the $j^{th}$ circle.

\section{Filtering section 6 of \cite{Doub}}

\noindent Next we associate maps in the knot Floer homology to the changes in resolutions at crossings. In our case, these maps become maps between filtered groups. We work backwards from the maps in \cite{Doub}. \\
\ \\
\noindent First, we note that the resolution changes occur in three ways: between circles split from the axis, between circles linking the axis, and
between circles of mixed linking. The first occur precisely as in \cite{Doub} due to the local nature of the surgeries in the double cover and the
connected sum decomposition of the covering manifolds. In particular, the maps for the filtered theory are just the maps for the unfiltered theory tensored with the identity on the tensor products of the $V$-vector spaces. Hence they reflect the differential of the reduced Khovanov homology. \\
\ \\
\noindent Now consider a resolution change joining two circles which link the axis. In the double cover, this corresponds to a cobordism which 
involves $0$-surgery on a curve which is homologically non-trivial and intersects only those spheres intersecting the binding. Such a cicle is isotopic to a circle in a fiber of the open book determined by $\# B(0,0)$ before connect summing with extra $S^{1} \times S^{2}$'s. Moreover, since the circle is the lift of an arc between two branch points, it is homologically non-trivial in the fiber. Ignoring the choice of basis implicit in the above description, we can calculate the effect of such a surgery by looking at the standard picture of $B(0,0)$ and doing $0$ surgery on a meridian of one of the $0$-surgered components of the Borromean rings. When we connect sum with copies of $B(0,0)$ we obtain a diffeomorphic picture to the one described above. We then use homology classes to pin down the maps in the original picture. In the unfiltered version, the model calculation
uses the following long exact sequence (which must split as depicted due to ranks and gradings). 

$$
\begindc{\commdiag}[10]
\obj(0,0)[I1]{$\cdots$}
\obj(20,0)[Ad1]{$\mathbb{F}_{\frac{1}{2}} \oplus \mathbb{F}_{-\frac{1}{2}}$}
\obj(10,0)[Ad]{$\mathbb{F}_{1} \oplus \mathbb{F}^{2}_{0} \oplus \mathbb{F}_{-1}$}
\obj(30,0)[Ad2]{$\mathbb{F}_{\frac{1}{2}} \oplus \mathbb{F}_{-\frac{1}{2}}$}
\obj(40,0)[I2]{$\cdots$}
\mor{I1}{Ad}{$-\frac{1}{2}$}
\mor{Ad}{Ad1}{$-\frac{1}{2}$}[1,5]
\mor{Ad1}{Ad2}{$0$}[1,1]
\mor{Ad2}{I2}{$- \frac{1}{2}$}
\enddc
$$

\noindent In the identification with Khovanov homology, the $\mathbb{F}_{1}$ term corresponds to $v_{+} \otimes v_{+}$ and it thus maps to $w_{+}$. The term mapping to $\mathbb{F}_{-\frac{1}{2}}$ in the surjection is the image under $\nu_{2}$ of $\mathbb{F}_{1}$ where $\nu_{2}$ is the meridian we \textit{do not} surger. Meanwhile the image of $\nu_{1}$ is annihilated. Transfering back to the basis from the resolution, this tells us that $\gamma_{1} + \gamma_{2}$ generates the kernel, and $\gamma_{1}$ and $\gamma_{2}$ are mapped isomorphically to $\gamma'$. Transfering back further to Khovanov's notation, we get $\gamma_{1} \ra v_{-} \otimes v_{+} \ra w_{-} \leftarrow \gamma'$ and $v_{+} \otimes v_{-} \ra w_{-}$.  \\
\ \\
\noindent For the filtered version, we obtain the model long exact sequence which filters the above one:

$$
\begindc{\commdiag}[5]
\obj(0,0)[I0]{$\ $}
\obj(0,8)[I1]{$\cdots$}
\obj(0,16)[I2]{$\ $}
\obj(10,0)[Ad]{$\mathbb{F}_{-1}$}
\obj(10,8)[Bd]{$\mathbb{F}_{0} \oplus \mathbb{F}_{0}$}
\obj(10,16)[Cd]{$\mathbb{F}_{1}$}
\obj(20,8)[Ad1]{$\mathbb{F}_{\frac{1}{2}} \oplus \mathbb{F}_{-\frac{1}{2}}$}
\obj(30,0)[Ad2]{$\mathbb{F}_{-\frac{1}{2}}$}
\obj(30,8)[Bd2]{$\mathbb{F}_{\frac{1}{2}} \oplus \mathbb{F}_{\frac{1}{2}}$}
\obj(30,16)[Cd2]{$\mathbb{F}_{\frac{3}{2}}$}
\obj(40,0)[T0]{$\ $}
\obj(40,8)[T1]{$\cdots$}
\obj(40,16)[T2]{$\ $}
\mor{I0}{Ad}{$\ $}
\cmor((12,7)(13,6)(17,6)(20,6)(21,7)) \pup(18,7){$\ $}
\cmor((18,9)(19,10)(23,10)(27,10)(28,9)) \pdown(28,12){$\ $}
\mor{I1}{Bd}{$\ $}
\mor{I2}{Cd}{$\ $}
\mor{Ad2}{T0}{$\ $}
\mor{Bd2}{T1}{$\ $}
\mor{Cd2}{T2}{$\ $}
\enddc
$$

\noindent The first term is the knot Floer homology of $B(0,0)$; the second is the knot floer homology of the unknot in $S^{1} \times S^{2}$ which we obtain after the $0$-surgery on the meridian; the third term is the result of $+1$-surgery on the meridian, $B(0,-1)$ in the notation of \cite{Knot}. The grading and ranks again determine the filtered maps on the first page. When we join two curves which link the axis, we obtain one which does not link the axis. This can be seen by considering the possible winding numbers for the result: 0 or 2. However, the result is a Jordan curve in the plane and thus cannot have winding number 2 about the origin. Working back through the basis transformations as before, these correspond in our notation to the maps $v_{+} \otimes v_{+} \ra 0$, $v_{+} \otimes v_{-} \ra w_{-}$, $v_{-} \otimes v_{+} \ra w_{-}$ and $v_{-} \otimes v_{-} \ra 0$. \\
\ \\
\noindent Finally, the model calculation in the cases of joining a linked with an unlinked circle corresponds to the map in the following diagram:

\begin{center}
\includegraphics[scale=0.9]{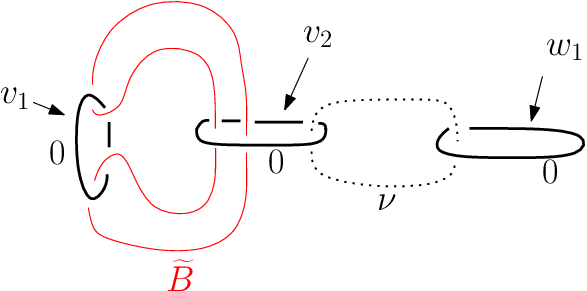}
\end{center}

\noindent The surgery circle, $\nu$, annihilates $\gamma_{1} + \gamma_{2}$ again in mod-2 homology. The result of the resolutioin change is now a circle which links the axis. The relevant cobordism map is from $B(0,0) \# S^{1} \times S^{2} \ra B(0,0)$ and corresponds to $v_{+} \otimes w_{+} \ra v_{+}$, $v_{+} \otimes w_{-} \ra 0$, $v_{-} \otimes w_{+} \ra v_{-}$ and $v_{-} \otimes w_{-} \ra 0$. This can be seen from the following graded exact sequence:
$$
\begindc{\commdiag}[6]
\obj(0,0)[I0]{$\ $}
\obj(0,8)[I1]{$\cdots$}
\obj(0,16)[I2]{$\ $}
\obj(10,0)[Ad]{$\mathbb{F}_{-\frac{3}{2}} \oplus \mathbb{F}_{-\frac{1}{2}}$}
\obj(10,8)[Bd]{$\mathbb{F}_{-\frac{1}{2}}^{2} \oplus \mathbb{F}_{\frac{1}{2}}^{2}$}
\obj(10,16)[Cd]{$\mathbb{F}_{\frac{1}{2}}\oplus \mathbb{F}_{\frac{3}{2}}$}
\obj(20,0)[Ad1]{$\mathbb{F}_{-1}$}
\obj(20,8)[Bd1]{$\mathbb{F}_{0}^{2}$}
\obj(20,16)[Cd1]{$\mathbb{F}_{1}$}
\obj(30,0)[Ad2]{$\mathbb{F}_{-1}$}
\obj(30,8)[Bd2]{$\mathbb{F}_{0}^{2}$}
\obj(30,16)[Cd2]{$\mathbb{F}_{1}$}
\obj(40,0)[T0]{$\ $}
\obj(40,8)[T1]{$\cdots$}
\obj(40,16)[T2]{$\ $}
\mor{I0}{Ad}{$\ $}
\mor{I1}{Bd}{$\ $}
\mor{I2}{Cd}{$\ $}
\mor{Ad}{Ad1}{$\ $}[1,5]
\mor{Bd}{Bd1}{$\ $}[1,5]
\mor{Cd}{Cd1}{$\ $}[1,5]
\mor{Ad2}{T0}{$\ $}
\mor{Bd2}{T1}{$\ $}
\mor{Cd2}{T2}{$\ $}
\enddc
$$
where $\mathbb{F}_{\frac{3}{2}}$ corresponds to $v^{1}_{+} \otimes v^{2}_{+} \otimes w_{+}$ and is mapped to $v^{1}_{+} \otimes v'_{+}$, taking into account
both $0$-framed knots in the Borromean rings. Note that a $w_{-}$ always forces the map to be $0$.  
\\
\ \\
\noindent Due to the introduction of the two new components we do not need to examine what happens if one of the circles is the marked circle: a division or merging never includes the marked circle.   \\
\ \\
\noindent Similar considerations, or duality, establish the maps for the case of splitting a circle into two circles. Note that the above maps are from $+1$ resolutions to $0$ resolutions. This force us to use the mirror of $\mathbb{L}$ in establishing the relationship between the knot Floer homology of $\widetilde{B}$ and the reduced skein homology. 

\begin{prop}
Let $\mathcal{P}$ be a projection for $\mir{\mathbb{L}}' \cup B$. Let $R$ be a choice of resolution for each crossing of $\mir{\mathbb{L}}'$. Then there is an isomorphism 
$$
\Phi_{B}(R) : \widetilde{V}(\mathcal{P}(R), B) \stackrel{\cong}{\lra} \widehat{HFK}(\Sigma(\mathcal{P}(R)), \widetilde{B})
$$
Let $R'$ be a resolution found by changing a single smoothing in $R$ from $0$ to $+1$. Then the following diagram commutes
$$
\begindc{\commdiag}[6]
\obj(0,0)[A1]{$\widehat{HFK}(\Sigma(\mathcal{P}(R)), \widetilde{B})$}
\obj(20,0)[A2]{$\widehat{HFK}(\Sigma(\mathcal{P}(R')), \widetilde{B})$}
\obj(0,10)[A3]{$\widetilde{V}(\mathcal{P}(R), B) $}
\obj(20,10)[A4]{$\widetilde{V}(\mathcal{P}(R'), B) $}
\mor{A3}{A1}{$\Phi_{B}(R)$}
\mor{A4}{A2}{$\Phi_{B}(R')$}
\mor{A1}{A2}{$\widehat{F}_{R<R'}$}
\mor{A3}{A4}{$d_{\mir{\mathbb{L}}}$}
\enddc
$$
\ \\
\noindent where $\widehat{F}_{R<R'}$ is the cobordism map for the knot Floer homologies induced by the surgery corresponding to the resolution change, and
$d_{\mir{\mathbb{L}}}$ is the differential in the skein homology. This square is a $\Z$-direct sum of squares where the index corresponds to the filtration of the knot Floer homology and the $k$ index in the skein homology (with $\mathcal{F} = \frac{k}{2}$ after the final shifts).
\end{prop}

\section{The knot Floer homology spectral sequence}

\noindent Let $L = L_{1} \cup \ldots \cup L_{n}$ be a framed link in a three manifold $Y$. Following section 4 of \cite{Doub} we let $R = (m_{1}, \ldots, m_{n})$ where $m_{i} \in \{0, 1, \infty\}$ and $Y(R)$ be the result of $fr(L_{i}) + m_{i}\,\mu_{i}$ surgery on each $L_{i}$ where $\mu_{i}$ is the meridian of $L_{i}$ and $\infty$-surgery is $\mu_{i}$-surgery. We let $0 < 1 < \infty$ define a lexicographic ordering on $\{0,1,\infty\}^{n}$ and call $I'$ and immediate successor to $I$ as in \cite{Doub}: all the $m_{j}'$ are the same as $m_{j}$ except for one where $m_{i}' > m_{i}$, excluding
the case $(m_{i}',m_{i}) = (\infty, 0)$. Then to each immediate successor $I'$ of $I$ there is a map
$$
F_{R<R'} : \widehat{HF}(Y) \lra \widehat{HF}(Y')
$$
arising from the associated two handle addition. \\
\ \\
\noindent According to section 8 of \cite{Knot}, two handle additions attached in a way algebraically unlinked from a knot induce maps on knot Floer homologies. Viewed differently, the knot turns the chain map above into a filtered morphism for the filtered homology groups. The ``top'' levels of these filtered morphisms form exact sequences which specialize to the skein exact sequence. Following these thoughts leads to
 
\begin{prop}
Let $L = L_{1} \cup \ldots \cup L_{n}$ be a framed link in $(Y, K)$ such that $\mathrm{lk}(L_{s}, K) =0$ for all $s$. For each integer $k$, and surface $F$ spanning $K$ and disjoint from $L$, there is a spectral sequence such that
\ben
\item The $E^{1}$ page is $\oplus_{R \in \{0,1\}^{n}} \widehat{HFK}(Y(R),K, k)$
\item The $d_{1}$ differential is obtained by adding all $\widehat{F}_{R<R'}$ where $R'$ is an immediate successor of $R$ 
\item All the higher differentials respect the dictionary ordering of $\{0,1\}^{n}$, and
\item The spectral sequence eventually collapses to a group isomorphic to $\widehat{HFK}(Y, K, k)$. 
\een
\end{prop}
\ \\
\noindent {\bf Proof:} Combine section 8 of \cite{Knot} with section 4 of \cite{Doub} to construct a complex $X$ which records the maps induced from the surgeries described above. Use the approach in section 8 to adjust the maps in section 4 to reflect the additional basepoint. Use lemma 8.1 of \cite{Knot} to ensure that the maps are filtration preserving. It is noted in \cite{Knot} that the various chain homotopies are also filtration preserving. This allows us to use the lemma found in the appendix, a filtered version of proposition 4.4 in \cite{Doub}. Now check the model calculations in section 4 of \cite{Doub} to see that they still cancel appropriately, a result of the linking number vanishing. The argument in Theorem 4.7 establishes the result when $L$ consists of a single component. For the induction in the proof of Theorem 4.1 of \cite{Doub} we note that the adjusted statement is that $E^{1}(X)$ and $E^{1}(X(S))$ are isomorphic to $0$, by applying the stated reasons to the $E^{1}$ maps. This allows us to repeat the induction and conclude that the the knot Floer complex is filtered quasi-isomorphic to an iterated mapping cone, namely the analog of the surgery complex built in \cite{Doub}, which will be denoted $X(\mathbb{L})$. $\Diamond$\\
\ \\
\noindent We consider the link in a connect sum of $S^{1} \times S^{2}$'s upon which we perform the surgeries corresponding to the resolution changes.
These are all algebraically split from the binding, so the above proposition applies for the choice of framing coming from the crossing data. Combining
with the proposition from the previous section, and adjusting for the filtration information, allows us to identify the $E^{1}$-page with the skein homology complex. As in \cite{Doub} we may then conclude  

\begin{prop}
Let $\mathbb{L}$ be a link in $A \times I \subset \R^{2} \times \R$ as above. There is a spectral sequence whose $E^{2}$ term is isomorphic to the reduced Khovanov skein chain complex of $\overline{\mathbb{L}}'$ in $A \times I$ with coefficients in $\mathbb{F}_{2}$ and which converges to $\oplus_{k \in \Z} \widehat{HFK}(\Sigma(\mathbb{L}) \#^{2} (S^{1} \times S^{2}), \widetilde{B} \# B(0,0), k , \mathbb{F}_{2})$. 
\end{prop}

\noindent By splitting according to the filtration data we can obtain the slightly stronger

\begin{prop}
There is a spectral sequence whose $E^{2}$ term is isomorphic to the sub-complex of the reduced Khovanov skein complex of $\overline{\mathbb{L}}'$ generated by the enhanced states with $\Psi(S) = 2k$ and which converges to $\widehat{HFK}(\Sigma(\mathbb{L}) \#^{2} (S^{1} \times S^{2}), \widetilde{B} \# B(0,0), k)$. 
\end{prop}

\noindent In fact, by taking the direct sum of all the groups for the knot Floer homologies over all the resolutions we can obtain a 
bi-filtered complex, filtered by the pair $(I, \Psi)$, where the $E^{1}$ term corresponds to the filtration of the bi-filtered reduced Khovanov
homology complex. Using the graded objects for just the $\Psi$ filtration and taking their homology produces the first proposition above. The additional terms in the maps in the Khovanov complex induce maps in the $E^{2}$ level of the spectral sequence using the $\Psi$ filtration, since
these correspond to terms in the filtered cobordism maps between the Heegaard-Floer homologies. These maps fit together to provide a filtered version of the spectral sequence in section 4 of \cite{Doub} with $K$ inducing the filtration. Additional pages ultimately calculate the Heegaard-Floer homology of the branched double cover.  \\
\ \\
\noindent More can be concluded from the proof outlined above and the homological algebra in the appendix.

\begin{lemma}\label{lem:spec}
For each $r \geq 1$, the $E^{r}$ page of the spectral sequence for $\widehat{HF}(\Sigma(\mathbb{L})\#^{2} (S^{1} \times S^{2}))$ computed from 
$\oplus_{k \in \Z} \widehat{HFK}(\Sigma(\mathbb{L}) \#^{2} (S^{1} \times S^{2}), \widetilde{B} \# B(0,0), k , \mathbb{F}_{2})$
using the differential from the knot Floer homology is quasi-isomorphic to the $E^{r}$ page for the filtered complex $X(\mathbb{L})$, computed using the maps induced from the link surgeries spectral sequence above.
\end{lemma}

\noindent Thus, the values of $\Psi$  will filter the Heegaard-Floer homology groups of the double branched cover in a way corresponding to
that induced by the spectral sequence for $\widetilde{B}$ (namely, the associated graded groups will be isomorphic). 

\section{Transverse links, open books and contact invariants}\label{sec:contact}

\noindent First, we note that 

\begin{thm} 
Any transverse link is transversely isotopic to a braid closure. Furthermore, two braids represent transversally isotopic links  if an only if one can be obtained from the other by  conjugations in the braid group, positive Markov moves, and their inverses.
\end{thm} 

\noindent This is the culmination of work by Bennequin for the first part, and by V. Ginzburg, S. Orevkov, and N. Wrinkle, who independently proved the second part. We will replace the contact structure with an open book. The standard contact structure on $S^{3}$ is supported by the open book with unknotted binding and discs for pages.  In the braid picture, this corresponds to including the axis of the braid, which is an unknot. When we take a branched cover of a transverse link, the contact structure lifts to a contact structure in the cover where we use a Martinet contact neighborhood of the transverse link. In the open book picture, this contact structure is supported by the pre-image of the open book, whose fibers are now more complicated, but whose binding is the lift of the axis. This follows since the lifted contact structure remains $C^{0}$-close to the pages of the open book, and transverse to the binding. We call this contact structure $\xi$. The contact structure on $\#^{2} \big(S^{1} \times S^{2}\big)$ induced by
the fibered knot $B(0,0)$ will be denoted $\xi_{0}$. \\
\ \\
\noindent For a braid, O. Plamenevskaya, \cite{Plam}, \cite{Plan}, defines a cycle, $\widetilde{\psi}(\mathbb{L})$, in the reduced Khovanov homology chain group. First she resolves all the crossings in the direction of the oriented braid. This constructs the maximal number of non-trivial loops in the skein algebra perspective. She then labels every one of the unmarked strands with a $-$ and the marked strand with a $+$. This enhanced state is closed in the reduced Khovanov homology theory, \cite{Plan}. \\
\ \\
\noindent Let $\mathbb{L}$ be a braid whose closure is the transverse link.  

\begin{theorem}
Suppose $\mathbb{L}$ intersects the spanning disc for $B$ an odd number of times. Then the element $\widetilde{\psi}(\mathbb{L}')$ is closed in the skein Khovanov homology and represents the unique homology class with minimal $\Psi$-grading. Under the correspondence with the $E^{2}$ term of the spectral sequence converging to knot Floer homology, it maps to an element which survives the spectral sequence and generates $\widehat{HFK}(-\Si(\mathbb{L}) \#^{2} \big(S^{1} \times S^{2}\big), \widetilde{B} \# B(0,0), -1-g(\widetilde{B})) \cong \mathbb{F}_{2}$. Upon mapping this last group into $\widehat{HF}(-\Si(\mathbb{L}) \#^{2} \big(S^{1} \times S^{2}\big))$, $\widetilde{\psi}(\mathbb{L}')$ corresponds to the contact element $c(\xi \# \xi_{0})$.
\end{theorem}

\noindent {\bf Note:} The correspondence at the end is not the same as first mapping $\psi$ into the reduced Khovanov homology and then considering the
spectral sequence from it to the Heegaard-Floer homology of the branched double cover. \\
\ \\
\noindent {\bf Proof:} There is only one element in the skein chain group which has $\Psi = -2g(\widetilde{B} \# B(0,0))$ and that is Plamenevskaya's element. For if $b$ is the braid index of $\mathbb{L}'$ then Euler characteristic calculations imply that $1 -2g(\widetilde{B} \# B(0,0)) = 2 - b$  and thus $\Psi$ must equal $1 - b$. This can only happen when all the crossings are resolved in the direction of the link so that there are $b$ non-trivial circles and precisely one circle (the marked one) is adorned with a $+$ sign. $\widetilde{\psi}(\mathbb{L}')$ is characterized as the unique enhanced state with minimal value for $\Psi(S)$ and thus generates the homology in this $f$-grading. This enhanced state survives in the spectral sequence for the knot Floer homology of the binding and yields in the limit a generator of $\widehat{HFK}(-\Sigma(\mathbb{L}'), \widetilde{B} \# B(0,0), -g(\widetilde{B} \# B(0,0))) \cong \mathbb{F}$ since it is the only generator in the filtration level. \\ 
\ \\
\noindent  The branched cover of $B$ over $\mathbb{L}'$ is $\widetilde{B} \# B(0,0)$ which supports the contact structure $\xi \# \xi_{0}$.  The contact element $c(\xi \# \xi_{0})$ is the image in $\widehat{HF}(-\Sigma(\mathbb{L}'))$ of the generator of $\widehat{HFK}(-\Sigma(\mathbb{L}'), \widetilde{B} \# B(0,0), -g(\widetilde{B} \# B(0,0)))$. Lemma \ref{lem:spec} guarantees that this generator corresponds to the 
$-g(\widetilde{B} \# B(0,0))$ level of the associated graded group for $\widehat{HF}(-\Sigma(\mathbb{L}'))$. This level is either $\cong \mathbb{F}$
or $\cong 0$ depending upon whether the contact element vanishes. Thus, Plamanevskaya's element converges to the contact element in the Heegaard-Floer
homology (with $\mathbb{F}_{2}$ coefficients).  $\Diamond$ \\

\begin{cor}
Under the correspondence in the previous theorem, $\widetilde{\psi}(\mathbb{L})$ corresponds to $c(\xi) \in \widehat{HF}(-\Si(\mathbb{L}), \mathbb{F}_{2})$.  
\end{cor}
\ \\
\noindent {\bf Proof:} If $\mathbb{L}$  intersects the spanning disc for $B$ an even number of times, use a positive Markov move to increase the
number of strands by 1. $\widetilde{\psi}(\mathbb{L})$ is mapped to $\widetilde{\psi}(\mathbb{L}_{+})$, \cite{Plan}, under this move. Meanwhile, in the double cover this corresponds to positively stabilizing the open book, and thus does not change the contact invariant. Renaming $\mathbb{L}_{+}$ by $\mathbb{L}$ we may now assume $\mathbb{L}$ intersects the spanning disc an odd number of times. Furthermore, $\widetilde{\psi}(\mathbb{L})$ clearly corresponds to $\widetilde{\psi}(\mathbb{L}')$ in a precise way. Using the previous theorem we have that $\widetilde{\psi}(\mathbb{L}')$ maps to $c(\xi \#\xi_{0})$ under the spectral sequence and $c(\xi \# \xi_{0}) = c(\xi) \otimes c(\xi_{0})$. Adding the two meridional strands tensors both homologies with $V^{\otimes 2}$. Thus, $c(\xi)$ in the knot Floer homology of $\widetilde{B}$ corresponds to Plamenevskaya's element in the skein homology of $\mathbb{L}$ since both are alterred in the same formal manner by the introduction of the new strands. $\Diamond$\\

\noindent We now turn to proving a the non-vanishing result mentioned in the introduction. We begin with a lemma:

\begin{lemma}
Let $\mathcal{C}$ be a bifiltered complex over a field. Then up to isomorphism there is a unique
bifiltered complex $\mathcal{C}'$ such that
\ben
\item $\mathcal{C}'$ is bifiltered chain homotopy equivalent to $\mathcal{C}$
\item $\mathcal{C}'_{ij} \cong H_{\ast}(\mathcal{C}_{ij})$
\item The differential $d'=\sum d'_{ij}$ on $\mathcal{C}'$ has $d_{00}' = 0$, and induces
the same spectral sequences for both filtrations. 
\een 
\end{lemma}

\noindent {\bf Proof:} Use the cancellation lemma as per sections 4 and 5 of Rasmussen's thesis, but only for
those elements with the same bifiltration indices. $\Diamond$\\
\ \\
\noindent We note that since the knot Floer spectral sequence for $\#^{k} B(0,0)$ collapses at $E^{2}$, the
use of the above lemma for the $I$-filtration means that $\oplus_{j} \mathcal{C}_{ij}'$ is isomorphic to
the knot Floer homology for the summands in the cube complex corresponding to that $I$-value. In particular, 
there are no differentials keeping $I$ fixed, and reducing $\Psi$. For lack of a better name, we will also call this reduced complex $X(\mathbb{L})$, or just $X$.
As a result, $E^{1}_{I}(X) \cong X$ for the filtration from $I$. Since $X$
is bi-filtered chain homotopy equivalent to $X(S)$, it too is quasi-isomorphic to the chain complex
for $\widehat{CF}(-\Si(\mathbb{L}))$ by a $\Psi$-filtered map.   \\
\ \\
\noindent We begin with a little notation: we let $X_{j}$ be the sub-complex of $X$ with $\Psi \leq j$. Likewise, let
$K_{j}$ be the sub-complex of the reduced Khovanov homology with the same condition. Now the $I$-filtration -- from the 
flattened cube-- filters these sub-complexes and their quotient complexes.

\begin{cor}
Suppose there exists a $n$ such that
\ben
\item $\psi(\mathbb{L})$ is exact in $K_{n}$
\item The $I$-induced spectral sequence on $X_{n}/X_{-2g}$ collapses at $E^{2}$.
\een
then $c(\xi) = 0$.
\end{cor}

\noindent The second condition, of course, makes some complex computed from the knot Floer chain groups isomorphic to the
corresponding complex computed from the skein Khovanov chain groups\footnote{In an earlier version of this paper, the author incorrectly asserted that the vanishing of $\psi$ is enough to conclude that $c(\xi)$ also vanishes. John Baldwin, \cite{Bald}, pointed out the error and has since discovered examples where $c(\xi)$ is non-zero despite $\psi$ vanishing in the reduced Khovanov homology.}\\
\ \\
\noindent {\bf Proof:} Suppose $\psi(\mathbb{L})$ has the bifiltration value $(I_\psi, -2g)$. If we try to compute
the homology of $X_{n}$ using the $I$-filtration, then $\psi$ generates the only group 
in the $\Psi$-filtration level $-2g$. Since $\psi$ is exact in $K_{n}$, there is some element with $I$-filtration 
$I_{\psi} - 1$ whose differential in $K_{n}$ is $\psi$ (recall the differential {\em increases} I-values). This element, $\nu$, may be a linear combination of
elements with many different $\Psi$ values. We note that $\nu$ is closed and not exact in $E^{1}(X_{n}/X_{-2g})$ as a chain complex computing $E^{2}$. It is closed since the only non-zero portion of $\partial_{Kh}\nu$ is in $X_{-2g}$. It is not exact since it would need to be the differential of something with higher $I$-filtration, and for those elements the differential, which is given by the Khovanov differantial, is the same as in $E^{1}(X_{n})$; however in $E^{1}(X_{n})$, $\nu$ is not closed and hence is not exact. Thus $[\nu]$ will be non-zero in $E^{2}(X_{n}/X_{-2g})$. \\
\ \\
\noindent Consider $C_{i}$ to be the sub-complex of $X_{n}$ with $I$-filtration greater than or equal to $i$. We have 
the commutative diagram represented in Figure \ref{fig:CommD}, to which the remaining argument refers. 
Here $\F$ is the homology of $X_{-2g}$, $Q_{c}$ is the quotient complex of $X_{n}$ by $C_{I_{\psi}+1}$ and $Q$ is the quotient complex
by $\{\Psi \leq -g\} \cup \{I \geq I_{\psi}+1\}$. The $0$ in the upper left comes from the observation that there are no generators in $X$ with
$\Psi \leq -2g$ and $I \geq I_{\psi}+1$. The $0$ on the map in the upper right indicates that it will generate the trivial map in homology 
due to $\nu$. From now on we let $X' = X_{n}/X_{-2g}$. 

\begin{figure}
$$
\begindc{\commdiag}[6]
	\obj(0,10)[I01]{$0$}
	\obj(0,20)[I02]{$0$}
	\obj(10,0)[I10]{$0$}
	\obj(10,10)[I11]{$C_{I_{\psi}+1}$}
	\obj(10,20)[I12]{$C_{I_{\psi}+1}$}
	\obj(10,30)[I13]{$0$}
	\obj(20,0)[I20]{$0$}
	\obj(20,10)[I21]{$X_{n}/X_{-2g}$}
	\obj(20,20)[I22]{$X_{n}$}
	\obj(20,30)[I23]{$\F$}
	\obj(20,40)[I24]{$0$}
	\obj(30,0)[I30]{$0$}
	\obj(30,10)[I31]{$Q$}
	\obj(30,20)[I32]{$Q_{C}$}
	\obj(30,30)[I33]{$\F$}
	\obj(30,40)[I34]{$0$}
	\obj(40,10)[I41]{$0$}
	\obj(40,20)[I42]{$0$}
	\obj(40,30)[I43]{$0$}
	\mor{I01}{I11}{$ $}
	\mor{I02}{I12}{$ $}
	\mor{I13}{I12}{$ $}
	\mor{I12}{I11}{$\mathrm{Id}$}
	\mor{I11}{I10}{$ $}
	\mor{I13}{I23}{$ $}
	\mor{I12}{I22}{$ $}
	\mor{I11}{I21}{$ $}
	\mor{I24}{I23}{$ $}
	\mor{I23}{I22}{$ $}
	\mor{I22}{I21}{$ $}
	\mor{I21}{I20}{$ $}
	\mor{I23}{I33}{$\mathrm{Id}$}
	\mor{I22}{I32}{$ $}
	\mor{I21}{I31}{$ $}
	\mor{I34}{I33}{$ $}
	\mor{I33}{I32}{$0$}[1,1]
	\mor{I32}{I31}{$ $}
	\mor{I31}{I30}{$ $}
	\mor{I33}{I43}{$ $}
	\mor{I32}{I42}{$ $}
	\mor{I31}{I41}{$ $}
\enddc
$$ 
\caption{\ }\label{fig:CommD}
\end{figure}
\ \\
\noindent An element of $H_{\ast}(Q)$ in filtration level $I_{\psi}-1$ must
have non-trivial representative in $E^{2}(X_{n}/X_{-2g})$. Furthermore, the argument above shows that $[\nu] \neq 0$
in $H_{\ast}(Q)$. This is certainly true in $Q_{c}$ since $\nu$ has a non-trivial differential. However, if in $Q$ there
is an element with differential equal to $\nu$, the only other possibility is that in $Q_{c}$ this element has
differential equal to $\nu$ plus something in $X_{-2g}$. But then $\partial^{2} \neq 0$ on this element. \\
\ \\
\noindent Suppose, $[\nu]$ has non-zero image, $[\omega]$, under the map $H_{\ast}(Q) \ra H_{\ast}(C_{I_{\psi}+1})$.
If $[\omega]$ has non-zero image in $H_{\ast}(X_{n})$, from the middle row, then it too must have a non-zero representative
in $E^{2}(X_{n}/X_{-2g})$, since $C_{I_\psi+1}$ has no representatives with $\Psi$-filtration $-2g$. But then the induced
differential from the long exact sequence implies that $\partial[\nu]=[\omega]$ in $X'$. Hence, the rank of $H_{\ast}(X')$
is strictly less than that at $E^{2}$, i.e. there is a non-trivial differential beyond $E^{2}$.\\
\ \\
\noindent Given the assumptions, we must have $[\omega] = 0$ in $H_{\ast}(X_{n})$. Then it is the image of some non-zero element, $[\eta]$ 
of $H_{\ast}(Q_{c})$. This element injects into $H_{\ast}(Q)$ so that $[\nu]$ minus the image of $[\eta]$ is the image 
of some non-zero element of $H_{\ast}(X')$. Furthermore, since the image of $[\nu]$ under $H_{\ast}(Q) \ra \F$ is non-zero, then
the map $H_{\ast}(X') \ra \F$ is surjective. As a result, the map $\F \ra H_{\ast}(X_{n})$ is zero,
but this implies that $H_{\ast}(X_{-2g}) \ra H_{\ast}(X)$ is zero.\\
\ \\
\noindent The filtered quasi-isomorphism from $X$ to $\widehat{CFK}$ induces a commutative diagram
$$
\begindc{\commdiag}[8]
\obj(0,10)[A]{$H_{\ast}(X_{-2g}) \cong \F$}
\obj(15,10)[B]{$H_{\ast}(X)$}
\obj(0,0)[C]{$HFK(\widetilde{B}, -g) \cong \F$}
\obj(15,0)[D]{$HF(-\Si(\mathbb{L}))$}
\mor{A}{B}{$\cdot 0$}
\mor{C}{D}{$ $}
\mor{A}{C}{$\cong$}
\mor{B}{D}{$\cong$}
\enddc
$$
thereby showing that the contact invariant vanishes. 
$\Diamond$\\
\ \\
\noindent {\bf Note \#1:} The purpose of the $\mathbb{F}_{2}$ coefficients is to connect with the extant versions of Khovanov homology.  In the end, the crucial observation is that Plamenevskaya's element uniquely defines the lowest filtered portion of the both the skein and knot Floer homologies. As long as this remains true  and there is an analogous reduced Khovanov homology, the same argument will work with other coefficients. In particular, just changing the sign conventions will not change the conclusion, but there should be some sign convention lifted from the Heegaard-Floer world which will allow $\Z$-coefficients. \\
\ \\
\noindent {\bf Note \#2:} For a braid, $\mathbb{L}$, we can lift a negative crossing or a positive crossing to negative/positive Dehn twists along 
homologically non-trivial curves in the fiber of the open book. These, in turn, fit into the long exact sequences of Heegaard-Floer and knot Floer 
homology. One sign fits into the $\infty, 0, +1$ sequence for the fiber framing, while the other fits into the $-1, 0, \infty$ sequence. Doing all the
surgeries at the same time yields a spectral sequence as in the previous section, with the maps in the $E^{1}$ page coming either from the maps
$\infty \ra 0$ or $0 \ra \infty$ from the respective long exact sequence and converging to the appropriate homology of the fibered knot. This is
the same sequence as that constructed above, only the basis for the framings has been alterred. Namely, if the framing from the crossing is declared
$\infty$ and the crossing is negative, then the $0$ framing is $\infty$ in the fiber framing, and $+1$ is $0$ in the fiber framing. The knot for the surgery is the same, a lift of an arc between two branched points.  \\
\ \\
\noindent We now collect some results for quasi-positive braids.  We note that for a quasi-positive braid, the lifted contact structure is Stein fillable. We can use the above argument to reprove that the induced contact element is non-vanishing, \cite{Cont}. Let the braid be given by $w_{1}\sigma_{i_{1}}w_{1}^{-1} \cdots$ $w_{k}\sigma_{i_{k}}w_{k}^{-1}$. We resolve only those crossings corresponding to the $\sigma_{i_{k}}$ terms. For the $00\ldots0$ resolution, the result will be $b$ non-trivial circles. Any $1$ resolutions make the situation more difficult, but all the non-zero terms occur in higher filtration levels. Plamenevskaya's element is then in the lowest level of the $00\ldots0$ resolution. There is no possibility in the spectral sequences of a higher differential landing at this spot as they must all map to enhanced states with at least one $1$ in their code. Thus the element survives in this spectral sequence. We now note that when there is a $1$ in the code, and the resulting resolution does not consist of unlinked circles, that the Heegaard-Floer homology of its double cover is, as a filtered group, the limit of a spectral sequence. Combining all of these shows that Plamenevskaya's element survives in the spectral sequence and thus gives the non-triviality of the contact element in the double branched cover. 
 
\section{Knot Floer results for alternating branch loci}

\noindent We now explore the implications of the previous sections for knot Floer homology. First, we define some notation. Let $\mathbb{L}$ be a link in $A \times I$ admitting a connected, alternating projection to $A$. According to \cite{Doub}, the Heegaard-Floer homology of $\Si(\mathbb{L}, \SP{s})$ is congruent to $\mathbb{F}$ for each of the $Spin^{c}$ structures on $\Si(\mathbb{L})$. For a $Spin^{c}$ structure $\SP{s}$ and a null-homologous knot $K \subset \Si(\mathbb{L})$ define

$$
\tau(K, \SP{s}) = \min_{s \in \Z} \big\{ s : \widehat{HF}(\mathcal{F}_{s}, \SP{s}) \stackrel{i_{\ast}}{\lra} \widehat{HF}(\Si(\mathbb{L}), \SP{s})\mathrm{\ is\ nontrivial}\big\}
$$   
\ \\
\noindent where $\mathcal{F}_{s}$ is the sub-complex of generators with filtration index less than or equal to $s$. 
Using the results of the previous sections and Lemma \ref{lem:spec} we can prove

\begin{theorem}\label{thm:alex}
Let $\mathbb{L}$ be a non-split alternating link in $A \times I$ intersecting the spanning disc for $B$ in an odd number of points. Then for each $k$ there is an isomorphism
$$
\widehat{HFK}(-\Si(\mathbb{L}) \#^{2} \big( S^{1} \times S^{2} \big), \widetilde{B} \# B(0,0), k) \cong \bigoplus_{i,j \in \Z} H^{i;j, 2k}(\mathbb{L})
$$
where, for each $Spin^{c}$ structure, the elements on the right side all have the same absolute $\Z/2\Z$-grading. Together these isomorphisms induce a filtered quasi-isomorphism from the $E^{2}$-page of the knot Floer homology spectral sequence to that of the skein homology spectral sequence. Thus the knot Floer spectral sequence collapses after two steps. Furthermore, for any $\SP{s} \in Spin^{c}(\Si(\mathbb{L}))$ we have that
$$
\tau(\widetilde{B}, \SP{s})  = 0
$$
where $\widetilde{B}$ is considered in $\Si(\mathbb{L})$.
\end{theorem}

\noindent The content of this theorem is that all the knots $\widetilde{B}$ have the same knot Floer properties as alternating knots in $S^{3}$,
and their knot Floer homology (over all $Spin^{c}$ structures) is determined by the skein homology, and the G$\ddot{\mathrm{o}}$retz matrix of $\mathbb{L}$, when applicable. The last is used to calculate the signature, through a formula of C. Gordon and R. Litherland, and the Heegaard-Floer invariants, $d(\SP{s})$, for $\SP{s} \in Spin^{c}(-\Si(\mathbb{L}))$, \cite{Doub}, which determine the precise absolute grading for the homology groups. However, it seems difficult to recover data about individual $Spin^{c}$ structures from the Khovanov formalism. \\
\ \\
\noindent {\bf Proof of theorem \ref{thm:alex}:} We have established that there is a spectral sequence starting at the right side of the isomorphism and converging
to the left side. The right side is the $E^{2}$ page of this spectral sequence. The $E^{3}$ page is computed using maps between resolutions
differing in at least two positions. Thus the maps will necessarily increase the $i$ grading by $2$. However, the right side is supported in
those triples satisfying $k - j + 2i = \sigma(\mathbb{L})$. If $i$ increases by $\delta_i$ and $k$ stays fixed, then $j$ must increase by $2\delta_i$. We compare
this to the grading in the knot Floer cubical spectral sequence which reflects the absolute grading shifts in the long exact sequences. For example,
in the Heegaard-Floer long exact sequence for surgery on an unknot, we have $\Z_{0} \ra \Z_{-\frac{1}{2}} \oplus \Z_{\frac{1}{2}}$. Thinking
of this as arising from a $0 \ra 1$ resolution change, we would shift the right side so that the $q$-gradings would be preserved. In our case,
we would shift up by $\frac{1}{2}$. In general, $q$ corresponds to a shift by $\frac{1}{2}$ in the absolute gradings. In the mapping cone construction, we shift the right side down by $1$ so that the chain map contributes to a $-1$ differential. Due to all our three manifolds being connect sums of $S^{1} \times S^{2}$'s and our surgery circles as either creating or destroying one of the summands, this same calculation can be applied throughout the cubical complex to obtain a coherent relative grading. Keeping track of the shifts yields that we can measure the relative grading in Heegaard-Floer homology by $\frac{\Delta\,j}{2} - \Delta\,i$ where $\Delta$ is a change in
the specified index. For this grading the differential at the $E^{1}$ page is a $-1$-differential. More importantly, once we are at the $E^{2}$ page, if we fix $k$ and consider triples with $k - j + 2i = \sigma(\mathbb{L})$ we see that the relative grading is $0$ between generators on this plane. All higher differentials must be $-1$-differentials, so the spectral sequence collapses in each $k$ grading at $E^{2}$. After the $E^{2}$ page if $\frac{\Delta\,j}{2} - \Delta\,i = -1$ then $\Delta\,k = -2$. Note that this is the shift in the additional terms defining the spectral sequence 
converging to Khovanov homology. \\
\ \\
\noindent This is not the difference in the absolute gradings on $-\Si(\mathbb{L})$ as the $Spin^{c}$ structures have different invariants, $d(\SP{s})$.
However, it does return the difference between $\widetilde{gr}({\bf x}) - d(\SP{s}) - 1$ for generators corresponding to $\SP{s}({\bf x}) = \SP{s}$. To do this measure from the element which generates $\widehat{HF}(-\Si(\mathbb{L}), \SP{s})$ tensored with $\Theta^{++} \in \widehat{HF} ( \#^{2} \big(S^{1} \times S^{2}\big))$. The first must exist in a single $(i,j)$-pair by the results of E. S. Lee, \cite{ELee}, and the comments in the next paragraph. Note that the generator of $\widehat{HF}(-\Si(\mathbb{L}), \SP{s})$ lies in the even absolute $\Z/2\Z$-grading, and thus the relative
grading above will be correct for Euler characteristic calculations, up to sign. \\
\ \\
\noindent The spectral sequence on the reduced skein homology collapses at its $E^{2}$-term (the $E^{3}$ term in the sequence we
are considering in this proof). At that stage we recover $\widehat{HF}(-\Si(\mathbb{L}))$, as the reduced Khovanov homology has total rank
given by $\mathrm{det}(\mathbb{L})$. By Lemma \ref{lem:spec}, the spectral sequence on the Khovanov skein homology is quasi-isomorphic to 
that on the knot Floer homology of $\widetilde{B} \# B(0,0)$ in $-\Si(\mathbb{L})$. This allows us to draw the conclusion concerning $\tau$. 
Namely, the Heegaard-Floer homology of $-\Si(\mathbb{L}) \#^{2} \big(S^{1} \times S^{2}\big)$ will have the form 
$\widetilde{H} \otimes V^{\otimes 2}$ and will lie on four lines $j - 2i = -\sigma(\mathbb{L})$ (with multiplicity 2) and $j - 2i = - \sigma(\mathbb{L}) \pm 2$. When we factor out the $V^{\otimes 2}$, we have the reduced homology lying on $k = \sigma(\mathbb{L}) + j - 2i = 0$. 
Since there is only one grading in each filtration level in the knot Floer homology, this implies that $\tau(\mir{\widetilde{B}}) = 0$ 
from which the result follows. $\Diamond$\\
\ \\
\noindent We can also derive some information about $\widetilde{B}$ for the branch loci depicted in Figure \ref{fig:special}, regardless of whether $\mathbb{L}$ is alternating. In particular, $\widetilde{B} \subset S^{3}$ in these cases and

\begin{lemma}
For $\widetilde{B}$ coming from the branch loci depicted in Figure \ref{fig:special},
$$
\tau(\widetilde{B}) = - \frac{1}{2} T(\mathbb{L})
$$
\end{lemma}

\noindent {\bf Proof:} Add two non-trivial, non-interacting unknots to $\mathbb{L}$ and mark one of these. There is then a spectral sequence
converging to the knot Floer homology of $\widetilde{B} \#^{2} B(0,0)$ from $H^{\ast;\ast\ast}(\mir{\mathbb{L}}) \otimes V$. Consider an element in the subcomplex corresponding to knot filtrations less than or equal to  $\tau(\widetilde{B}) - 1$ which maps to $\Theta^{--} \in \widehat{HF}(\#^{2} S^{1} \times S^{2})$ under inclusion of the subcomplex. Then there is a element with $k$-gradings less than or equal to $2\tau(\widetilde{B}) - 2$ which survives the spectral sequence to the knot Floer homology . However, since $\mathbb{L}$ is an unknot, $\Theta^{--}$ is the element
 ${\bf u}_{-1} \otimes v_{-}$. Therefore, this same element will survive the spectral sequence from the skein homology to the Khovanov homology. Hence
 $T(\mir{\mathbb{L}}) - 2 \leq 2\tau(\widetilde{B}) - 2$ and $- \frac{1}{2} T(\mathbb{L}) \leq \tau(\widetilde{B})$. This is also true for $\mir{B}$ whence $- \frac{1}{2}T(\mir{\mathbb{L}}) \leq \tau(\mir{\widetilde{B}})$. Therefore, $- \frac{1}{2}T(\mathbb{L}) \geq \tau(\widetilde{B})$ as well. $\Diamond$\\
\ \\
\noindent These results hold in slightly greater generality. In the sequel to this paper an argument is given which holds for a broader class of links, similar to the quasi-alternating links of \cite{Doub}. This is the {\it smallest} subset of links in $A \times I$, denoted $\mathcal{Q}'$, with the property that
\ben
\item The alternating, twisted unknots, linking $B$ an odd number of times, are in $\mathcal{Q}'$.
\item 
 If $L \subset A \times I$ is a link admitting a connected projection to $A$, with a crossing such that
 \ben
 \item[] $\bullet$ The two resolutions of this crossing, $L_{0}$ and $L_{1}$, are in $\mathcal{Q}'$ and are
 connected in $A$, and
 \item[] $\bullet$ $\mathrm{det}(L) = \mathrm{det}(L_{0}) + \mathrm{det}(L_{1})$
 \een
 then $L$ is in $\mathcal{Q}'$
\een
The alternating $L$ used above are in $\mathcal{Q}'$, and the elements of $\mathcal{Q}'$ when considered in $S^{3}$ are all quasi-alternating as in \cite{Doub}. For this class of links Wehrli's algorithm terminates at the base cases of our induction, from which the conclusion in the theorem can be drawn. For braids in $\mathcal{Q}'$ we can be more precise about Plamenevskaya's element:\\

\begin{cor}
Let $\mathbb{L}$ be in $\mathcal{Q}'$. If the element $\widetilde{\psi}$ vanishes in the reduced Khovanov homology, then $c(\xi) = 0$.
\end{cor}
\ \\
\noindent {\bf Proof:} This corollary follows from the non-vanishing result in Section 8 since the spectral sequence for $X_{2-2g}/X_{-2g}$
collapses according to theorem \ref{thm:alex}. However we need to verify that, $\widetilde{\psi}$ is zero
in $Kh_{1-g}$. The only difficulty arises if there is a $\nu$ whose Khovanov differential is $\widetilde{\psi}$ and $\nu = \sum \nu_{i}$ where
$\nu_{i}$ is in $\Psi$-filtration level $2i$. Since the Khovanov differential reduces $\Psi$ by at most $2$, this requires the $i$ indices to range
from $1-g$ to $l$, and for there to be a summand for each index in the range. As we collapse the complex along differentials preserving the $\Psi$-filtration level, the complex stabilizes at $E^{2}$, and the structure described above yields a differential from $\nu_{l}$ to $\widetilde{\psi}$. However, we know that at $E^{2}$ $k - j + 2i = \sigma(\mathbb{L})$, and $\nu_{l}$ and $\widetilde{\psi}$ must have the same
$j$ value since they are linked by Khovanov differentials. In addition, the change in $i$ is an increase of $1$ from $\nu_{l}$ to $\widetilde{\psi}$. This implies that $k$ must decrease by $2$, and thus $l = 1-g$ as required.  $\Diamond$  \\  

\section{Examples}

\begin{figure}
\begin{center}
\includegraphics[scale=0.6]{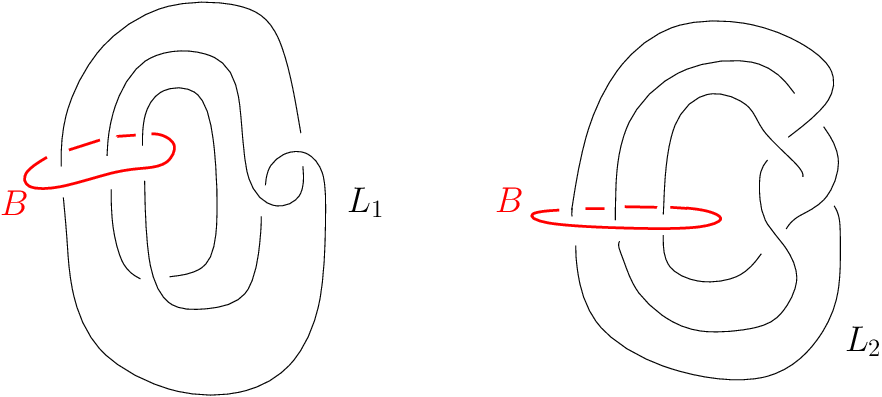}
\end{center}
\caption{The diagram for example 1 is on the left; that for example 2 is on the right.}\label{fig:examp1}
\end{figure}

\noindent {\bf Example 0:} Let $\mathbb{L}$ be a non-split alternating link and suppose $B$ is a meridian of one of the components. Then
$\widetilde{B}$ is an unknot in $-\Si(\mathbb{L})$ since the spanning disc lifts to a disc. Mark the link as above, then the reduced skein homology
after the final shifting agrees with the reduced Khovanov homology. On the other hand, the knot Floer homology of this unknot is just
$\widehat{HF}(-\Si(\mathbb{L}))$ in filtration level $0$. The equivalence of these two groups is a consequence of \cite{Doub}. In this sense, 
theorem \ref{thm:alex} is a generalization of the result in \cite{Doub}.\\
\ \\
\noindent {\bf Example 1:} See Figure \ref{fig:examp1} for the diagram. Here $\mathbb{L}$ is an unknot in $S^{3}$, so $\widetilde{B}$ is a knot
in $S^{3}$ as well. Untwisting and taking the branched double cover (or using symmetry between the two components) shows that $\widetilde{B}$ is 
the knot:

\begin{center}
\includegraphics[scale=0.6]{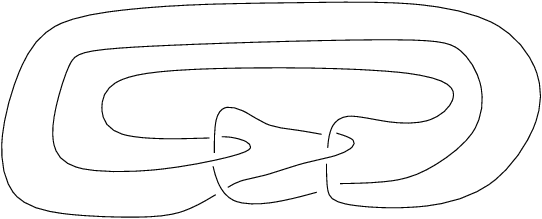}
\end{center}

\noindent This is the alternating knot, $6_{1}$, with signature equal to $0$. The main result in \cite{Alte} now verifies  the knot Floer conclusions of theorem \ref{thm:alex}. Furthermore, the Alexander polynomial is $-2\,T^{-1} + 5 - 2\,T$. We content ourselves with a direct verification of the
rank of the highest filtration level. Only resolutions with three non-trivial circles contribute to this level. These resolutions and the associated generators are: \\

\begin{center}
\includegraphics[scale=0.6]{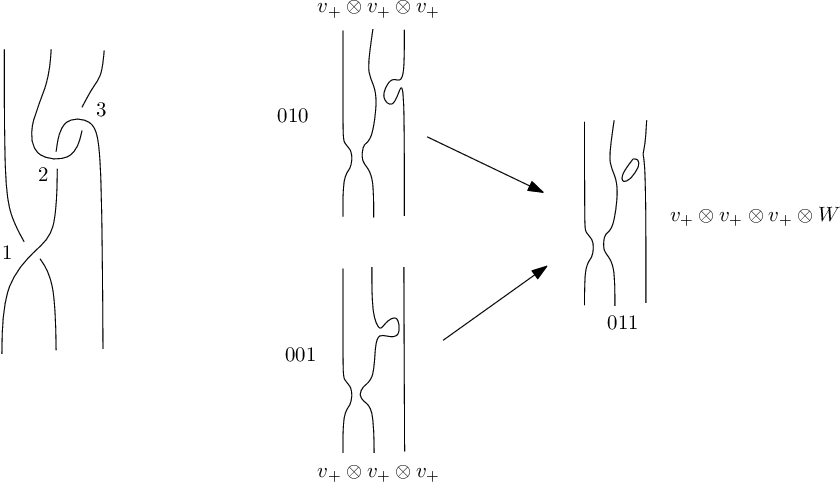}
\end{center}

\noindent The maps from $010$ and $001$ to $011$ both take $v_{+}\otimes v_{+} \otimes v_{+}$ to $v_{+} \otimes v_{+} \otimes v_{+} \otimes w_{-}$, and
thus their sum is closed, as is $v_{+} \otimes v_{+} \otimes v_{+} \otimes w_{+}$. The latter is two $q$-gradings above the former closed element, but
it also has one more $1$-resolution. Shifting $q$ down by $2$ decreases the homological grading by $1$ when identifying with knot Floer homology. Thus,
these generators are in the same grading in the knot Floer complex. This confirms theorem \ref{thm:alex} for the highest filtration level (modulo some
shifting).\\
\ \\
\noindent {\bf Example 2:} See Figure \ref{fig:examp1} for the diagram. Here $\mathbb{L}$ is the figure 8 knot, $4_{1}$, whose branched double cover is $L(5,2)$. In this arrangement, $\widetilde{B}$ is a genus $1$ fibered knot in $L(5,2)$. The possibilities for such a knot are strictly limited,
since there is only a $\Z$ in filtration levels $\pm 1$. The real content of the theorem here is that $\tau(\widetilde{B}) = 0$, as this
implies that there is one $Spin^{c}$ structure where the knot Floer homology is that of $4_{1}$. We give a non genus $1$ example later.\\
\ \\
\noindent The monodromy for this knot is $\big(\gamma_{1}\gamma_{2}^{-1}\big)^{2}$ where
$\gamma_{i}$ is a positive Dehn twist around a standard symplectic basis element for $H_{1}(T^{2} - D^{2})$. The monodromy action on $H_{1}$ and the
Alexander polynomial associated to the $\Z$-covering from the fibering are computed to be 
$$
A = \left[\begin{array}{cc} 2 & 3 \\ 3 & 5 \\ \end{array} \right]  \hspace{.5in} \Longrightarrow \hspace{.5in} \mathrm{det}(I - tA) \doteq  \Delta_{\widetilde{B}}(t)  = - T^{-1} + 7 - T^{1}
$$
where we have symmetrized and normalized $\mathrm{det}(I - tA)$ according to the convention in the second appendix. We can compute the skein homology directly, but instead we use the theory to compute it from the polynomial $V(t, q, x)$. This is not quite as direct as it may seem. The polynomial satisfies
$$
V(t, q, x)/ (qx + q^{-1}x^{-1}) = t^{-2}q^{-4} + t^{2}q^{4} + 2\,t^{-1}q^{-2} + 2\,tq^{2} + q^{2}x^{2} + 1 + q^{-2}x^{-2} + t^{1} + 1
$$
From our conventions, we should add two non-trivial strands, and at the end factor out $V^{\otimes 2}$ to get to the knot Floer homology. However,
adding a marked non-trivial circle and a non-trivial circle amounts to multiplying $qx(qx + q^{-1}x^{-1})V(t, q, x)$, so will
equal the above after shifting and removing $V^{\otimes 2}$.\\
\ \\
\noindent Note that the second to last term does not satisfy $k - j + 2i = 0$ and so should not appear in the homology. It cancels with one of the $1$'s. In fact, this can be seen by resolving crossings in turn, until one arrives at a the closure of $\sigma_{2} \in B_{3}$, where we use the homology calculation from previously (the same cancellation occurs there) and then build up the homology independently. When all is completed
we obtain the following diagram on the $(j, k)$-axes (the subscripts are the values of $i$).
$$
\begindc{\commdiag}[4]
\obj(-5,5)[A2]{$-3$}
\obj(-5,10)[A3]{$-2$}
\obj(-5,15)[A4]{$-1$}
\obj(-5,20)[A5]{$0$}
\obj(-5,25)[A6]{$+1$}
\obj(-5,30)[A7]{$+2$}
\obj(-5,35)[A8]{$+3$}
\obj(-5,40)[A9]{$+4$}
\obj(5,0)[b2]{$-3$}
\obj(10,0)[b3]{$-2$}
\obj(15,0)[b4]{$-1$}
\obj(20,0)[b5]{$0$}
\obj(25,0)[b6]{$+1$}
\obj(30,0)[b7]{$+2$}
\obj(35,0)[b8]{$+3$}
\obj(40,0)[b9]{$+4$}
\obj(0,20)[b2]{$\mathbb{F}_{-2}$}
\obj(10,20)[b2]{$\mathbb{F}_{-1}^{2}$}
\obj(20,20)[b2]{$\mathbb{F}_{0}$}
\obj(30,20)[b2]{$\mathbb{F}_{1}^{2}$}
\obj(40,20)[b2]{$\mathbb{F}_{2}$}
\obj(10,30)[b2]{$\mathbb{F}_{-2}$}
\obj(30,10)[b2]{$\mathbb{F}_{2}$}
\enddc
$$
Note that if we shift the elements to $j =0$, decreasing $i$ by $1$ each time $j$ decreases by $2$, every group in the same horizontal row shifts to  the same grading. Note also that the ranks after shifting horizontally reflect the coefficients of the Alexander polynomial; and, up to a minus sign, the $\Z/2\Z$-gradings are correct. Furthermore, if we consider the $\Psi$-filtration, in the $E^{\infty}$-page of the spectral sequence there will be five terms on the $k = 0$ horizontal line, correponding to the five $Spin^{c}$ structures on $L(5,2)$. All that remains is to identify which generators correpond to which $Spin^{c}$ structure and then use the Goeritz matrix for $4_{1}$ to complete the absolute grading calculations. To do
this we should use the more refined torsion, $\check{\tau}(Y - K)$, in our Euler characteristic computations, \cite{Grig}. We complete this argument
in the sequel to this paper. Comparing the two will show that the $\Z/2\Z$-gradings from the knot Floer homology 
correspond to those from the skein homology. However, the correspondence only occurs when we add over all $Spin^{c}$ structures and all $q$-gradings.\\

\newpage

\appendix 

\section{Homological Algebra}

\noindent All coefficients are taken in $\mathbb{F}_{2}$, hence the difference from the usual signs. However, everything can be adapted to
work with coefficients in $\Z$.\\
\ \\
\noindent Let $(A,\mathcal{A})$ and $(B,\mathcal{B})$ be filtered differential modules. Let $f: A \ra B$ be a filtered chain map. Then the mapping
cone $M(f)$ inherits a filtration by declaring $\mathcal{M}_{i} = \mathcal{A}_{i} \oplus \mathcal{B}_{i}$. That the differential preserves this
filtration follows from $f$ being filtered. When undeclared, a filtration on a mapping cone complex will come from this construction. The definitions imply that $E^{1}(\mathcal{M}) \cong \mathrm{MC}(E^{1}(f))$.\\
\ \\
\noindent A filtered chain map $f$ will be a 1-quasi-isomorphism if it induces an isomorphism between the $E^{1}$ pages of the spectral sequences for the source and the target. For the morphism of spectral sequences induced by $f$, in which the induced maps intertwine the differentials on each page, this implies that all the higher pages, $E^{r}$, are quasi-isomorphic by the induced map, $E^{r}(f)$. This is probably weaker than
$f$ being a filtered chain isomorphism, but enough for spectral sequence computations.\\
\ \\
\noindent Let $\left\{(A_{i}, \mathcal{A}_{i})\right\}_{i=0}^{\infty}$ be a set of filtered chain complexes with each filtration $\mathcal{A}_{i}$
being bounded and ascending:
$$
\mathcal{A}_{i} : \{0\} = A_{i}^{n_{i}} \subset \cdots \subset A^{j}_{i} \subset A^{j+1}_{i} \subset \cdots \subset A^{N_{i}}_{i} \cong A_{i}
$$
Let $\left\{ f_{i} : A_{i} \ra A_{i+1}\right\}$ be a set of chain maps satisfying:
\ben
\item $f_{i}$ is a filtered map for each $i$.
\item $f_{i+1} \circ f_{i}$ is filtered chain homotopic to $0$, i.e. there is a filtered map $H_{i}: A_{i} \ra A_{i+2}$ such that 
$f_{i+1} \circ f_{i} = \partial_{i+2} \circ H_{i} + H_{i} \circ \partial_{i}$. 
\item $f_{i + 2} \circ H_{i} + H_{i+1}\circ f_{i} : A_{i} \ra A_{i+3}$ is an 1-quasi-isomorphism. 
\een
\ \\
\noindent In this setting we have the lemma, following \cite{Doub},

\begin{lemma} 
The mapping cone $\mathrm{MC}(f_{2})$ is 1-quasi-isomorphic to $A_{4}$.
\end{lemma}

\noindent {\bf Proof:} The hypotheses above guarantee that the maps in the proof of lemma 4.4 of \cite{Doub} are filtered maps.  We need only check the filtering condition for maps in and out of the mapping cone, but with the aforementioned convention these are clearly filtered. In particular
the map $\psi_{i} = f_{i + 2} \circ H_{i} + H_{i+1}\circ f_{i}$ is a 1-quasi-isomorphism by assumption, and the same argument as in \cite{Doub}
implies that $\alpha_{2}$ is a quasi-isomorphism which is also filtered. This is not quite enough to conclude, but it does ensure that $\alpha_{i}$ induces maps at each page in the spectral sequence.\\
\ \\
\noindent The module $Gr(A_{i}) \cong \oplus_{j\in \Z} A_{i}^{j}/A_{i}^{j-1}$ inherits 
a differential which maps the $j^{th}$ graded component to itself, whose homology provides $E^{1}$. The maps $f_{i}$ induce chain maps between these complexes for each grading level. Indeed each of the maps $\psi_{i}, H_{i}, f_{i}$, etc., likewise induce such maps. Compositions such as $f_{i+1} \circ H_{i}$ induce maps on the graded components which are the same as the compositions for the maps induced from $f_{i+1}$ and $H_{i}$ separately. Thus for each $j$, we have the situation in the lemma in \cite{Doub} applied solely to the $j^{th}$ graded component. Applying the lemma in each grading guarantees that the map induced in that grading by $\alpha_{2}$ is a quasi-isomorphism, i.e. that the induced map on the $E^{1}$ page is an isomorphism of spectral sequences. Thus, $\alpha_{2}$ induces an isomorphism from the $E^{1}$ page for $A_{4}$ to $\mathrm{MC}(E^{1}(f_{2})) \cong E^{1}(\mathrm{MC}(f_{2}))$, which is the desired result. $\Diamond$\\
\ \\
\noindent As in \cite{Doub}, we can reinterpret this as a result on interated mapping cones. Let $M = \mathrm{MC}(f_{1}, f_{2}, f_{3})$ be the 
filtered chain complex on $A_{1} \oplus A_{2} \oplus A_{3}$, filtered by $A_{1}^{j} \oplus A_{2}^{j} \oplus A_{3}^{j}$, and equipped with the
differential 
$$
\left(
\begin{array}{ccc}
\partial_{1} & 0 & 0 \\
f_{1} & \partial_{2} & 0 \\
H_{1} & f_{2} & \partial_{3} \\
\end{array}
\right)
$$
\ \\
\noindent That this is a differential is a consequence of the assumptions made before the lemma. The lemma then implies that the induced spectral 
sequence on the iterated mapping cone collapses at the $E^{1}$ term. This follows according to the following diagram:

$$
\begindc{\commdiag}[10]
\obj(0,0)[A00]{$\ $}
\obj(10,0)[A10]{$0$}
\obj(20,0)[A20]{$0$}
\obj(30,0)[A30]{$0$}
\obj(40,0)[A40]{$\ $}
\obj(0,3)[A01]{$0$}
\obj(10,3)[A11]{$A_{3}^{j}/A_{3}^{j-1}$}
\obj(20,3)[A21]{$M^{j}/M^{j-1} $}
\obj(30,3)[A31]{$\mathrm{MC}^{j}(f_{1})/\mathrm{MC}^{j-1}(f_{1})$}
\obj(40,3)[A41]{$0$}
\obj(0,6)[A02]{$0$}
\obj(10,6)[A12]{$A_{3}^{j}$}
\obj(20,6)[A22]{$A_{1}^{j-1} \oplus A_{2}^{j} \oplus A_{3}^{j} $}
\obj(30,6)[A32]{$A_{1}^{j} \oplus A_{2}^{j} $}
\obj(40,6)[A42]{$0$}
\obj(0,9)[A03]{$0$}
\obj(10,9)[A13]{$A_{3}^{j-1}$}
\obj(20,9)[A23]{$A_{1}^{j-1} \oplus A_{2}^{j-1} \oplus A_{3}^{j-1}$}
\obj(30,9)[A33]{$A_{1}^{j-1} \oplus A_{2}^{j-1}$}
\obj(40,9)[A43]{$0$}
\obj(0,12)[A04]{$\ $}
\obj(10,12)[A14]{$0$}
\obj(20,12)[A24]{$0$}
\obj(30,12)[A34]{$0$}
\obj(40,12)[A44]{$\ $}
\mor{A11}{A10}{$\ $}
\mor{A21}{A20}{$\ $}
\mor{A31}{A30}{$\ $}
\mor{A12}{A11}{$\ $}
\mor{A22}{A21}{$\ $}
\mor{A32}{A31}{$\ $}
\mor{A13}{A12}{$\ $}
\mor{A23}{A22}{$\ $}
\mor{A33}{A32}{$\ $}
\mor{A14}{A13}{$\ $}
\mor{A24}{A23}{$\ $}
\mor{A34}{A33}{$\ $}
\mor{A01}{A11}{$\ $}
\mor{A02}{A12}{$\ $}
\mor{A03}{A13}{$\ $}
\mor{A11}{A21}{$\ $}
\mor{A12}{A22}{$\ $}
\mor{A13}{A23}{$\ $}
\mor{A21}{A31}{$\ $}
\mor{A22}{A32}{$\ $}
\mor{A23}{A33}{$\ $}
\mor{A31}{A41}{$\ $}
\mor{A32}{A42}{$\ $}
\mor{A33}{A43}{$\ $}
\enddc
$$
where the top two rows are exact, and all the columns are exact. The nine lemma now guarantees that the bottom row is exact, and each of the maps
is a chain map. In the long exact sequence from the bottom row, there is one map guaranteed to be an isomorphism by the lemma. Consequently, the 
groups in $E^{1}(M)$ are trivial.

\end{document}